\documentclass{amsart}
\usepackage{amsfonts}
\usepackage{amsthm}
\usepackage{amsmath}
\usepackage[mathscr]{eucal}
\numberwithin{equation}{section}
\usepackage{nicefrac}
\usepackage{pgfplots, pgfplotstable}
\usepackage[font=large]{subcaption}

\newcommand{\uproman}[1]{\uppercase\expandafter{\romannumeral#1}}

\numberwithin{equation}{section}

\theoremstyle{plain}
\newtheorem{Th}{Theorem}[section]
\newtheorem{Lemma}[Th]{Lemma}
\newtheorem{Cor}[Th]{Corollary}

\theoremstyle{definition}
\newtheorem{Def}[Th]{Definition}
\newtheorem{Rem}[Th]{Remark}
\newtheorem{Ex}[Th]{Example}

\newcommand{\scp}[3]{\langle #1,#2\rangle_{#3}}
\newcommand{\minim}[1]{{\footnotesize\textsc{v}_{#1}}}
\newcommand{\minimt}[2]{{\footnotesize\textsc{v}_{#1}(#2)}}
\newcommand{\Vh}{\textsc{V}_h}
\newcommand{\Vr}{\mathcal{V}_r}
\newcommand{\Vnull}{\mathscr{V}_0}
\newcommand{\Veins}{\mathscr{V}_1}
\newcommand{\IntNorm}[2]{\|#1\|_{\textrm{#2}^s(\Vnull,\Veins)}}
\newcommand{\IntNormFEM}[2]{\|#1\|_{#2^s}}
\newcommand{\IntNormRB}[2]{\|#1\|_{#2_r^s}}
\newcommand{\IntSpace}[1]{[\Vnull,\Veins]_{\textrm{#1}^s}}
\newcommand{\TheIntSpace}{[\Vnull,\Veins]_{s}}
\newcommand{\Vbold}{\mathbb{V}{(\Vnull,\Veins; y^\alpha)}}
\newcommand{\R}{\mathbb{R}}
\newcommand{\N}{\mathbb{N}}
\newcommand{\Op}{\mathcal{L}}
\newcommand{\InducedScp}[3]{\scp{#1}{#2}{\textrm{#3}^s(\Vnull,\Veins)}}
\newcommand{\InducedFEMScp}[3]{\scp{#1}{#2}{#3^s}}
\newcommand{\InducedRBScp}[3]{\scp{#1}{#2}{#3_r^s}}
\newcommand{\InducedOp}[1]{\Op_{\textsc{#1}^s(\Vnull,\Veins)}}
\newcommand{\InducedFEMOp}[1]{\Op_{#1^s}}
\newcommand{\InducedRBOp}[1]{\Op_{#1_r^s}}
\newcommand{\InducedFEMMatOp}[1]{L_{#1^s}}
\newcommand{\InducedRBMatOp}[1]{L_{#1_r^s}}
\newcommand{\Krm}{\mathrm{K}}

\newcommand{\Norm}[2]{\|#1\|_{#2}}
\newcommand{\CoefVec}[1]{\underline{#1}}
\newcommand{\CoefvecVr}[1]{\underline{\underline{#1}}}

\DeclareMathOperator{\Span}{span}

\DeclareMathOperator{\Identity}{I}
\DeclareMathOperator{\Trace}{tr}
\DeclareMathOperator{\DeltaAmpl}{dn}

\begin{document}
	\title{A Reduced Basis Method For Fractional Diffusion Operators \uproman{1}}

	\author{Tobias Danczul}
	\author{Joachim Sch\"oberl}
	\address{TU Wien \\ Institute for Analysis and Scientific Computing \\
		Wiedner Hauptstrasse 8-10, 1040 Wien, Austria} 
	\email{tobias.danczul@tuwien.ac.at}
	\email{joachim.schoeberl@tuwien.ac.at}

	

	\begin{abstract} 
	We propose and analyze new numerical methods to evaluate fractional norms and apply fractional powers of elliptic operators. By means of a reduced basis method, we project to a small dimensional subspace where explicit diagonalization via the eigensystem is feasible. The method relies on several independent evaluations of $(\Identity-t_i^2\Delta)^{-1}f$, which can be computed in parallel. We prove exponential convergence rates for the optimal choice of sampling points $t_i$, provided by the so-called \textit{Zolotar\"ev points}. Numerical experiments confirm the analysis and demonstrate the efficiency of our algorithm.
	\end{abstract}
	
	\maketitle
	
	\section{Introduction}
	Fractional powers of differential operators are a field of substantial interest in different branches of mathematics. Their augmented appearance in real world problems, such as ecology \cite{ApplPopulDyn}, finance \cite{ApplApplebaum}, image processing \cite{ApplImageProc}, material science \cite{ApplMaterialSc}, and porous media flow \cite{ApplPorousMediaFlow} has given rise to several approaches in order to understand and analyze problems of this kind. 
	
	Typically, direct computations rely on matrix approximations $L$ of the desired operator, whose $s^{th}$-power is computed subsequently. This procedure requires diagonalization of $L$, which amounts to a large number of time-consuming eigenvalue problems, making this approach inapplicable for general purposes. Adding to this difficulty, many practical scenarios demand numerical methods that allow efficient evaluations in $s$. In \cite{sMapAppl2}, the fractional exponent is determined in a way, such that the observed data matches the mathematical model. In \cite{sMapAppl1}, $s$ serves as control parameter to minimize a given cost functional. All these investigations suggest that one is interested in the entire family of solutions for $s\in(0,1)$ rather than one specific value of $s$. The demand for suitable methods that address these problems has substantially increased throughout the last years.
	
	Fractional powers of the Laplace operator appear to be of particular interest. Widely varying definitions of $(-\Delta)^s$ have emerged, e.g., as pseudo differential operator defined by the Fourier transform, by means of its involved eigensystem, as a singular integral operator, or as inverse of the Riesz potential operator. All these definitions turn out to be equivalent in $\R^n$, see \cite{EquivDef}. This result no longer holds as bounded domains are incorporated. A detailed excursion of its versatile definitions as well as the comparison of both existing and newly proposed numerical schemes is explained in \cite{FracLaplaceReview} and \cite{FracDiffReview2}.
	
	A difficulty that all fractional operators have in common is their nonlocal character. Caffarelli and Silvestre managed to avoid this inconvenience in \cite{HarmonicExt} by relating any fractional power of the Laplacian in $\R^n$ to a Dirichlet-to-Neumann map of an involved harmonic extension problem in $\R^n\times\R^+$, providing a local realization of $(-\Delta)^s$. Adaptions for bounded domains $\Omega$ have been conducted in \cite{ConConvEllipProb}, \cite{SqrRootLaplaceDirichlet}, and \cite{RadialExtremalSol}, yielding a boundary value problem on the semi-infinite cylinder $\mathcal{C} := \Omega\times\R^+$. Enhancements for a more general class of operators has been presented in \cite{HarmExtGeneralized}.	
	
	A large number of methods exploits the structure of harmonic extension techniques to approximate fractional differential operators and their inverses, see \cite{AinsworthSpectralExt}, \cite{CertifiedRB}, \cite{MelenkTensorFEM}, and \cite{ExtensionProblemToSqrRoot}. In \cite{NochettoOtarolaSalgado}, the solution of the aforementioned boundary value problem is computed on the truncated cylinder $\mathcal{C}_{\gamma} := \Omega\times[0,\gamma)$, with $\gamma > 0$ of moderate size, by standard finite element techniques, at the cost of one additional space dimension. Truncation can be justified by the fact that the solution decreases exponentially in the extended direction. Other strategies rely on block-wise low-rank approximation of the associated stiffness-matrix \cite{MelenkHMatrix}. Of particular interest are fractional elliptic operators in context of parabolic equations. Tackling problems of this kind has been a matter of several recent publications, e.g., \cite{FracCahnHilliard}, \cite{PasciakParabolicFrac}, \cite{MelenkRiederer}, and \cite{SchwabNS}.
		
	In this article, we interpret fractional operators as interpolation operators and make use of the $\textrm{K}$-method \cite{RealMethoOfInt} to provide attractive approximations for any arbitrary symmetric, uniformly elliptic operator $\Op$. This approach requires the knowledge of the map $t\mapsto \minimt{}{t} := (I + t^2\Op)^{-1}f$, whose smoothness in $t$ justifies the usage of standard reduced basis technology, see e.g., \cite{RBMRef2}, \cite{RBMRef3}, and \cite{RBMRef}. Providing an optimal choice of snapshots $t_i$, we pursue a proper orthogonal decomposition (POD) strategy, see e.g., \cite{KunischVolkwein}, to project a matrix approximation $L$ of the desired operator to a lower dimensional space, where its fractional power can be determined directly. The decoupled structure of the projection trivially admits an efficient implementation in parallel. Multigrid preconditioner can be utilized, whose convergence rates are bounded independently of the shift parameter $t_i$ and uniform mesh size $h$. The proposed method can be interpreted as model order reduction of the approach devised in \cite{NochettoOtarolaSalgado}, without requiring truncation of the domain. Among others, it provides accurate approximations for evaluations of both types $s\mapsto(-\Delta)^su$ and $u\mapsto(-\Delta)^su$ with considerably reduced computational expense. The arising operator incorporates a nonlinear dependency in $u$. This inconvenience is compensated by its analytically affirmed exponential convergence property, while, at the same time, the computational effort grows only logarithmically in the condition number. Emerging estimates rely on rational approximation of the bivariate function $(1+\lambda^2t^2)^{-1}$ over a suitable spectral interval of the discrete operator. Realization of the inverse operator and parabolic problems that involve fractional diffusion operators  is the matter of a consecutive paper.
	
	We emphasize that this is not the first investigation that recognizes the importance of $\minimt{}{t}$ in context of fractional operators. The approach developed in \cite{Pasciak} and later improved in \cite{SincQuadImproved} relies on the Dunford-Taylor integral representation of $(-\Delta)^{-s}$. A sinc quadrature scheme especially tailored for integrals of this type is presented and requires the evaluation of $\minimt{}{t}$ at the quadrature nodes. Exponential decay of the error in the number of nodes is shown. Further modifications are discussed in \cite{RBMBonito} where computations of the quadrature are accelerated by means of a reduced basis method. While the construction of the reduced basis differs from the one we pursue, the involved proof of convergence also relies on rational approximations and partially follows the outline of our analysis. A similar approach is proposed in \cite{RBMKatoa} in consideration of a different quadrature rule, where exponential convergence rates are observed numerically. A model order reduction that relies on the extension method is proposed in \cite{CertifiedRB}.
	
	It remains to be mentioned that we are not the first to relate fractional powers of differential operators to rational approximation. In \cite{PasciakBURA}, for instance, the so-called best uniform rational approximation (BURA) of $t^{-s}$, which was originally proposed in \cite{BURA}, is utilized as matrix function to approximate $(-\Delta)^{-s}$.
	
	The paper is organized as follows. In Section \ref{SecInt}, we introduce three different concepts of Hilbert space interpolation. They serve as abstract template to provide a setting that is suitable for the study of our problem. Equivalence of all three methods appears to be the main result of this section. A reduced basis strategy is applied in Section \ref{SecNormAppr}. Feasible choices of the reduced space and the efficient implementation of its arising reduced basis interpolation norms are taken into account. Having understood its underlying structure, we proceed, in Section \ref{SecOpAppr}, to deduce the induced fractional operator, providing the essential definition of this paper. Numerical analysis is performed in Section \ref{Sec:ConvAna}. The optimal choice of the reduced space is elaborated, yielding exponential decay in the error for both norm and operator approximation. The core of this paper is summarized in Theorem \ref{Th:Core}. Thereupon, in Section \ref{SecNumExp}, we conduct several numerical examples that illustrate the performance of our method. Eventually, in the appendix, we prove two technical results that are referred to within Section \ref{SecInt}. 
	
	\section{Notation and Preliminaries}\label{SecInt}
	In this section, we establish the notation and terminology we utilize throughout this paper and introduce several function spaces we address in the subsequent.
	
	Throughout what follows, let the \textit{induced norm} $\Norm{\cdot}{}$ of a Hilbert space $(\mathscr{V},\scp{\cdot}{\cdot}{})$ be defined as 
	\begin{align}\label{SecInt:InducedDef}
		\Norm{\cdot}{} := \sqrt{\scp{\cdot}{\cdot}{}}.
	\end{align}
	Conversely, given a Banach space $(\mathscr{V},\Norm{\cdot}{})$, such that $\Norm{\cdot}{}$ satisfies the parallelogram law, we define the \textit{induced scalar product of} $\Norm{\cdot}{}$ as the unique scalar product $\scp{\cdot}{\cdot}{}$ on $\mathscr{V}$ that satisfies \eqref{SecInt:InducedDef}, obtained by polarization identity. Whenever referring to a Banach space $(\mathscr{V},\Norm{\cdot}{})$ as Hilbert space, we mean that $\Norm{\cdot}{}$ induces a scalar product $\scp{\cdot}{\cdot}{}$, such that $(\mathscr{V},\scp{\cdot}{\cdot}{})$ is a Hilbert space.
	\subsection{Hilbert space interpolation}\label{Sec:IntSpaces}
	Throughout what follows, let $(\Vnull,\scp{\cdot}{\cdot}{0})$ and $(\Veins,\scp{\cdot}{\cdot}{1})$ denote two real Hilbert spaces, such that $\Veins\subseteq \Vnull$ is dense with compact embedding. It is well-known that there exists an orthonormal basis $(\varphi_k)_{k=1}^\infty$ of $\Vnull$ and a sequence of positive real numbers $(\lambda_k)_{k=1}^\infty$ with $\lambda_k\longrightarrow\infty$ as $k\longrightarrow\infty$, satisfying
	\begin{align*}
		\llap{$\forall w\in \Veins:$\quad} \scp{\varphi_k}{w}{1} = \lambda_k^2\scp{\varphi_k}{w}{0}
	\end{align*}
	for all $k\in\N$. Along with these premises, we introduce, based on \cite{IntSpaces}, \cite{Bramble}, and \cite{NonHomBVPandAppl}, the first of three space interpolation techniques. For each $s\in(0,1)$ we define an interpolation space between $\Vnull$ and $\Veins$ by
	\begin{align*}
		\IntSpace{H} := \{u\in \Vnull: \IntNorm{u}{H}<\infty\},
	\end{align*} 
	equipped with its \textit{Hilbert interpolation norm}
	\begin{align*}
		\|u\|_{\textrm{H}^s(\Vnull,\Veins)}^2 := \sum_{k = 1}^{\infty}\lambda_k^{2s}u_k^2, \rlap{\qquad$u_k := \scp{u}{\varphi_k}{0}.$}
	\end{align*}
	$\left(\IntSpace{H},\IntNorm{\cdot}{H}\right)$ incorporates a Hilbert space structure and satisfies 
	\begin{align*}
		\Veins\subseteq \IntSpace{H} \subseteq \Vnull.
	\end{align*}
	
	Another approach is provided by the real method of interpolation in terms of the $\textrm{K}$-functional. It was first published by Peetre \cite{RealMethoOfInt}, Lions and Magenes \cite{NonHomBVPandAppl} and also works for Banach spaces. Let $\Norm{\cdot}{0}$ and $\Norm{\cdot}{1}$ denote the induced norms on $\Vnull$ and $\Veins$, respectively. We define for all $t>0$ and $u\in \Vnull$ the \textit{$\textrm{K}$-functional} as
	\begin{align*}
		\Krm_{(\Vnull,\Veins)}(t;u) := \inf\limits_{v\in \Veins}\sqrt{\|u-v\|_0^2 + t^2\|v\|_1^2}
	\end{align*}
	to obtain the \textit{$\textrm{K}$-norm}
	\begin{align*}
		\IntNorm{u}{K}^2 := \int_{0}^{\infty}t^{-2s-1}\textrm{K}_{(\Vnull,\Veins)}^2(t;u)\,dt.
	\end{align*}
	Along with its inner product, obtained by parallelogram law, the norm induces a Hilbert space
	\begin{align*}
		\IntSpace{K} := \{u\in \Vnull: \IntNorm{u}{K}<\infty\},
	\end{align*} 
	which again turns out to be intermediate.
	
	Based on the work of Peetre and Lions, it has been shown that $[\Vnull,\Veins]_{\textsc{K}^s}$ can be characterized as space of trace,  see \cite{IntToSob&IntSpac} and \cite{Triebel}. The arising norm, which turns the trace space into a Banach space, is known to be equivalent to the $\textrm{K}$-norm. We affirm this observation by proving that these norms are not only equivalent but do also coincide up to a multiplicative constant. This result is well-known for some particular choices of $\Vnull$ and $\Veins$, see e.g., \cite[Proposition 2.1]{RadialExtremalSol}, but has not been recorded in its most general setting to the best of our knowledge. To make matters precise, we investigate some technical results.
	
	For each $s\in(0,1)$ let $\alpha := 1-2s$ henceforth. For all $i = 0,1$, we define the space $L_2(\R^+,\mathscr{V}_i;y^\alpha)$ of all measurable functions $v:\R^+\longrightarrow \mathscr{V}_i$, such that
	\begin{align*}
		\int_{\R^+}y^\alpha\Norm{v(y)}{i}^2\,dy < \infty,
	\end{align*}
	and further
	\begin{align*}
		H^1(\R^+,\Vnull; y^\alpha) := \{v\in L_2(\R^+,\Vnull; y^\alpha):v'\in L_2(\R^+,\Vnull; y^\alpha)\}.
	\end{align*}
	Thereupon, we introduce
	\begin{align*}
		\Vbold := H^1(\R^+,\Vnull; y^\alpha)\cap L_2(\R^+,\Veins; y^\alpha)
	\end{align*}
	and endow it with the norm
	\begin{align*}
		 \|v\|_{\Vbold}^2  := \int_{\R^+}y^\alpha\left(\Norm{v(y)}{1}^2 + \Norm{v'(y)}{0}^2\right)\,dy.
	\end{align*}
	This space is amenable to trace evaluations, as the following Theorem shows.
	\begin{Th}\label{Th:Trace}
		There exists a linear, surjective trace operator 
		\begin{alignat*}{2}
			\Trace: \Vbold &\longrightarrow\; &&\IntSpace{\textsc{H}},\\
			v(y) &\mapsto &&v(0),
		\end{alignat*}
		such that for all $v\in\Vbold$ there holds
		\begin{align}\label{SecInt:TraceIneq}
			\sqrt{d_s}\IntNorm{\Trace v}{\textsc{H}}\leq\|v\|_{\Vbold}.
		\end{align}
		By $d_s$ we refer to a positive constant whose value can be specified by means of the Gamma function,
		\begin{align*}
			d_s = 2^{1-2s}\frac{\Gamma(1-s)}{\Gamma(s)}.
		\end{align*}
	\end{Th}
	\begin{proof}
		See Appendix.
	\end{proof}
	Theorem \ref{Th:Trace} justifies the introduction of an interpolation space
	\begin{align}\label{SecInt:ExtSpace}
		\IntSpace{E} := \Trace(\Vbold),
	\end{align}
	endowed with the \textit{extension-norm}
	\begin{align}\label{SecInt:ExtNorm}
		\IntNorm{u}{E} :=  \inf\limits_{\substack{\mathscr{U}\in \Vbold\\ \Trace \mathscr{U} = u}}  \|\mathscr{U}\|_{\Vbold},
	\end{align}
	which is is well-defined by standard arguments from calculus of variation. Due to surjectivity, \eqref{SecInt:ExtSpace} coincides with $\IntSpace{H}$.
	
	By means of Euler-Lagrange formalism, we observe that the infimum in \eqref{SecInt:ExtNorm} is the unique solution of an involved variational formulation. Thereupon, we introduce the following definition.
	\begin{Def}
		The $\alpha$-harmonic extension $\mathscr{U}$ of $u\in\IntSpace{E}$ is defined as the unique solution of the variational formulation: Find $\mathscr{U}\in \Vbold$, such that for all $y\in\R^+$ and $\mathscr{W}\in\Vbold$ there holds
		\begin{align}\label{SecInt:HarmonicExtPDE}
			\begin{split}
				 \scp{y^\alpha\mathscr{U}(y)}{\mathscr{W}(y)}{1} - \frac{\partial}{\partial y}\left(y^\alpha\scp{\mathscr{U}'(y)}{\mathscr{W}(y)}{0}\right) = 0,\\
				\Trace \mathscr{U} = u.
			\end{split}
		\end{align}
	\end{Def}
	\begin{Lemma}\label{Lm:Minimizer}
		Let $\mathscr{U}$ denote the $\alpha$-harmonic extension of $u\in\IntSpace{\textsc{E}}$. Then there holds
		\begin{align*}
			\IntNorm{u}{\textsc{E}} = \Norm{\mathscr{U}}{\Vbold}.
		\end{align*}
	\end{Lemma}
	\begin{proof}
		Follows directly from the fact that \eqref{SecInt:HarmonicExtPDE} is the Euler-Lagrange equation of the minimization problem in \eqref{SecInt:ExtNorm}.
	\end{proof}
	As the following Theorem shows, all three interpolation methods coincide.
	\begin{Th}\label{Th:NormEquiv}
		Let $\Vnull, \Veins$ denote two Hilbert spaces, such that $\Veins\subseteq \Vnull$ is dense with compact embedding. Then there holds 
		\begin{align*}
			\IntNorm{\cdot}{\textsc{E}} = \sqrt{d_s}\IntNorm{\cdot}{\textsc{H}} = \sqrt{d_s}C_s\IntNorm{\cdot}{\textsc{K}},
		\end{align*}
		where $C_s:= \sqrt{\frac{2\sin(\pi s)}{\pi}}$.
	\end{Th}
	\begin{proof}
		See Appendix.
	\end{proof}
	Throughout what follows, we denote by $\TheIntSpace$ the unique interpolation space emerging from one and hence from all three interpolation methods. The norms $\IntNorm{\cdot}{H}$ and $\IntNorm{\cdot}{K}$ satisfy the parallelogram law on $\TheIntSpace$. By virtue of Theorem \ref{Th:NormEquiv}, this property also applies to $\IntNorm{\cdot}{E}$. We summarize these observations in the following corollary. 
	\begin{Cor}\label{Cor:ScpEquiv}
		All three interpolation norms, $\IntNorm{\cdot}{\textsc{E}}$, $\IntNorm{\cdot}{\textsc{H}}$, and $\IntNorm{\cdot}{\textsc{K}}$, induce a respective scalar product, $\InducedScp{\cdot}{\cdot}{\textsc{E}}$, $\InducedScp{\cdot}{\cdot}{\textsc{H}}$, and $\InducedScp{\cdot}{\cdot}{\textsc{K}}$, such that for all $v,w\in \TheIntSpace$
		\begin{align*}
			\InducedScp{v}{w}{\textsc{E}} = d_s\InducedScp{v}{w}{\textsc{H}} = d_sC_s^2\InducedScp{v}{w}{\textsc{K}}.
		\end{align*}
	\end{Cor}
	\begin{Def}\label{Def:InducedOperator}
		Let $\Norm{\cdot}{s}\in\{\IntNorm{\cdot}{E}, \IntNorm{\cdot}{H}, \IntNorm{\cdot}{K}\}$. By means of the Riesz-representation Theorem, we define the induced operator of $\|\cdot\|_s$ as the unique linear function $\Op_s:\TheIntSpace\longrightarrow\TheIntSpace$, such that
		\begin{align*}
			\llap{$\forall v\in\TheIntSpace:$\quad} \scp{v}{\Op_sw}{0} = \scp{v}{w}{s}
		\end{align*}
		for each $w\in\TheIntSpace$, where $\scp{\cdot}{\cdot}{s}$ refers to the induced scalar product of $\Norm{\cdot}{s}$.
	\end{Def}
	Corollary \ref{Cor:ScpEquiv} immediately reveals the following result.
	\begin{Cor}\label{Cor:OpEquiv}
		All three interpolation norms, $\IntNorm{\cdot}{\textsc{E}}$, $\IntNorm{\cdot}{\textsc{H}}$, and $\IntNorm{\cdot}{\textsc{K}}$, induce an operator, $\InducedOp{\textsc{E}}$, $\InducedOp{\textsc{H}}$, and $\InducedOp{\textsc{K}}$, such that for all $v,w\in \TheIntSpace$
		\begin{align*}
			\scp{v}{\InducedOp{\textsc{E}}w}{0} = d_s \scp{v}{\InducedOp{\textsc{H}}w}{0} = d_sC_s^2\scp{v}{\InducedOp{\textsc{K}}w}{0}.
		\end{align*}
	\end{Cor}
	\begin{Rem}
		Let $\Omega\subseteq\R^d$, $d\in\N$, be a bounded domain with Lipschitz boundary. Along with the choice $\Vnull = (L_2(\Omega),\Norm{\cdot}{L_2})$ and $\Veins = (H_0^1(\Omega),\Norm{\nabla\cdot}{L_2})$, $\InducedOp{H}$ coincides with the spectral fractional Laplacian subject to homogeneous Dirichlet boundary conditions and $\InducedOp{E}$ with the involved Dirichlet-to-Neumann map from Caffarelli and Silvestre.
	\end{Rem}
	\subsection{The finite element framework}
	Results from Section \ref{Sec:IntSpaces} are also valid in a discretized setting. Depending on two fixed Hilbert spaces $\Vnull$ and $\Veins$ which satisfy the premises from Section \ref{Sec:IntSpaces}, we denote by $\Vh\subseteq \Veins$ a conforming finite element space of dimension $N$ henceforth. Further, let $(b_{k})_{k=1}^N\subseteq\Vh$ denote an arbitrary basis of $\Vh$. By $M, A\in\R^{N\times N}$, we refer to the mass- and stiffness-matrix of $\Vh$, arising from finite element discretization in terms of
	\begin{align}\label{SecInt:DefMassStiffMat}
		M_{ji} = \scp{b_i}{b_j}{0},\qquad A_{ji} = \scp{b_i}{b_j}{1}.
	\end{align}
	Due to its finite dimensional nature, the spaces $(\Vh,\|\cdot\|_{0})$ and $(\Vh,\|\cdot\|_{1})$ satisfy the conditions from Section \ref{Sec:IntSpaces}, such that the \textit{discrete interpolation norms on} $\Vh$
	\begin{align}\label{SecInt:DiscIntNorms}
		\begin{split}
			\IntNormFEM{u_h}{E} := &\|u_h\|_{\textsc{E}^s((\Vh,\|\cdot\|_{0}),(\Vh,\|\cdot\|_{1}))},\quad \IntNormFEM{u_h}{H} := \|u_h\|_{\textrm{H}^s((\Vh,\|\cdot\|_{0}),(\Vh,\|\cdot\|_{1}))}, \\
			&\IntNormFEM{u_h}{K}^2 := \int_{0}^{\infty}t^{-2s-1}\underbrace{{\Krm}_{((\Vh,\|\cdot\|_{0}),(\Vh,\|\cdot\|_{1}))}^2(t;u_h)}_{=: K^2(t;u_h)}\,dt,
		\end{split}
	\end{align}
	are well-defined. The finite element space equipped with each of these norms is a Banach space, inducing both scalar product, $\InducedFEMScp{\cdot}{\cdot}{E}$, $\InducedFEMScp{\cdot}{\cdot}{H}$, $\InducedFEMScp{\cdot}{\cdot}{K}$, and operator, $\InducedFEMOp{E}$, $\InducedFEMOp{H}$, $\InducedFEMOp{K}$, respectively. The aim of this paper is to provide an accurate approximation of these operators with considerably reduced computational expense. By virtue of Corollary \ref{Cor:OpEquiv}, it suffices to address this problem in any of those three interpolation settings. Each of them comes with its own benefits and difficulties attached. In the following section, we exploit the advantages of all three strategies to derive a computationally beneficial norm approximation, such that the induced operator satisfies the desired properties.

	\section{Approximation of the interpolation norms}\label{SecNormAppr}
	The goal of this section is to devise an accurate approximation of the discrete interpolation norms, introduced in \eqref{SecInt:DiscIntNorms}, with downsized computational effort. For convenience, we neglect the subscript $h$ for all finite element functions $u_h\in\Vh$ and solely write $u$ henceforth. Furthermore, by $(\varphi_k,\lambda_k^2)_{k=1}^N\subseteq \Vh\times\R^+$ we refer to the $0$-orthonormal eigenpairs of $(\Vh,\|\cdot\|_{0})$ and $(\Vh,\|\cdot\|_{1})$ from now on, such that
	\begin{align}\label{SecNormAppr:Eigenpairs}
		\llap{$\forall w\in \Vh:\qquad$} \scp{\varphi_k}{w}{1} = \lambda_k^2\scp{\varphi_k}{w}{0}.
	\end{align}
	The eigenvalues are assumed to be in ascending order according to their value, such that $$0<\lambda_1^2\leq...\leq\lambda_N^2.$$
	
	\subsection{The reduced basis approach}\label{SecNormAppr:SubSec:RedBasisAppr}
	Utilizing standard reduced basis technology for one	dimensional parametric elliptic partial differential equations, we define for each $u\in\Vh$ its approximate interpolation norms as follows.
	\begin{Def}[Reduced basis interpolation norms]\label{Def:Method}
		For each $t\in\R^+$ we define $\minimt{N}{t}\in\Vh$ as the unique solution of
		\begin{align}\label{SecNormAppr:ShiftLapl}
			\scp{\minimt{N}{t}}{w}{0} + t^2\scp{\minimt{N}{t}}{w}{1} = \scp{u}{w}{0}
		\end{align}
		for all $w\in\Vh$. Given some real parameters $0 = t_0 < t_1 < ... < t_r$, specified in Section \ref{Sec:ConvAna}, we introduce the reduced space
		\begin{align}\label{Def:ReducedSpace}
			\Vr :=\Span\{\minimt{N}{t_0},...,\minimt{N}{t_r}\}.
		\end{align}
		The reduced basis interpolation norms on $\Vr$ are defined by either of the three equivalent definitions 
		\begin{subequations}
			\begin{alignat}{2}
				&\IntNormRB{u}{E} &&:= \|u\|_{\textsc{E}^s((\Vr,\|\cdot\|_{0}),(\Vr,\|\cdot\|_{1}))}, \label{SecNormAppr:RBExtNorm}\\ 
				&\IntNormRB{u}{H} &&:= \|u\|_{\textrm{H}^s((\Vr,\|\cdot\|_{0}),(\Vr,\|\cdot\|_{1}))}, \label{SecNormAppr:RBHilbertNorm}\\
				&\IntNormRB{u}{K}^2 &&:= \int_{0}^{\infty}t^{-2s-1}\underbrace{\Krm_{\left((\Vr,\|\cdot\|_{0}),(\Vr,\|\cdot\|_{1})\right)}^2(t;u)}_{=: K_r^2(t;u)}\,dt. \label{SecNormAppr:RBKNorm}
			\end{alignat}
		\end{subequations}
	\end{Def}
	\begin{Rem}
		The choice $t_0 = 0$ yields $\minimt{N}{t_0} = u$ and hence $u\in\Vr$, such that \eqref{SecNormAppr:RBExtNorm}-\eqref{SecNormAppr:RBKNorm} are well defined. Definition \eqref{Def:ReducedSpace} is motivated by means of the $\textrm{K}$-method. The variational problem \eqref{SecNormAppr:ShiftLapl}, which is uniquely solvable according to Lax-Milgram, appears to be the Euler-Lagrange equation of $K^2(t;u)$, such that $\minimt{N}{t}$ coincides with the minimizer of $K^2(t;u)$. Based on a sophisticated selection of $t_1,...,t_r$, the choice of $\Vr$ aims to provide a both accurate and efficient approximation to the family of solutions $(\minimt{N}{t})_{t\in\R^+}$.
	\end{Rem}
	\begin{Rem}
		Definition \ref{Def:Method} incorporates a nonlinear dependency in $u$.  For simplicity, we neglect this relation in both terminology and notation throughout our discussions. We point out, however, that all $\Vr$-connected constructions are subject to this dependency.
	\end{Rem}
	In analogy to \eqref{SecNormAppr:Eigenpairs}, we denote the eigenpairs of $(\Vr,\Norm{\cdot}{0})$ and $(\Vr,\Norm{\cdot}{1})$ by $(\phi_j,\mu_j^2)_{j=0}^r\subseteq \Vr\times\R^+$ from now on, such that 
	\begin{align}\label{SecNormAppr:RbEigenpairs}
		\scp{\phi_j}{\phi_i}{0} = \delta_{ji},\qquad \scp{\phi_j}{w_r}{1} = \mu_j^2\scp{\phi_j}{w_r}{0}, \rlap{\quad$w_r\in\Vr,$}
	\end{align}
	with
	\begin{align*}
		0<\mu_0^2\leq ...\leq \mu_r^2.
	\end{align*}
	
	In general, the construction of $\Vr$ yields a $r+1$ dimensional space. The proof is carried out in two steps. 
	\begin{Lemma}\label{Lm:InfRepr}
		For all $t\in\R^+$ there holds
		\begin{align*}
			 \minimt{N}{t} = \sum\limits_{k = 1}^{N}\frac{u_k}{1+t^2\lambda_k^2} \varphi_k,\rlap{\qquad$u_k := \scp{u}{\varphi_k}{0}.$}
		\end{align*}
	\end{Lemma}
	\begin{proof}
		Both $\minimt{N}{t}$ and $u$ provide expansions in the eigenbasis
		\begin{align*}
			\minimt{N}{t} = \sum\limits_{k=1}^{N}\scp{\minimt{N}{t}}{\varphi_k}{0}\varphi_k, \qquad u = \sum\limits_{k=1}^{N}\scp{u}{\varphi_k}{0}\varphi_k.
		\end{align*}
		Plugging into \eqref{SecNormAppr:ShiftLapl} with $w = \varphi_i$ yields
		\begin{align*}
			\scp{\minimt{N}{t}}{\varphi_i}{0} +t^2\lambda_i^2\scp{\minimt{N}{t}}{\varphi_i}{0} = \scp{u}{\varphi_i}{0} = u_i.
		\end{align*}
		Resolving the equation for $\scp{\minimt{N}{t}}{\varphi_i}{0}$ concludes the proof.
	\end{proof}
	\begin{Lemma}\label{Lm:LinInd}
		Let $|\{\lambda_k^2\in\{\lambda_1^2,...,\lambda_N^2\}:u_k\neq 0\}| =: m\in\N$ denote the number of pairwise distinct eigenvalues, whose corresponding eigenspaces contribute nontrivially to the linear combination of $u$. Then there holds
		\begin{align*}
			r+1\leq m \quad\Longrightarrow \quad\{\minimt{N}{t_0},..., \minimt{N}{t_r}\}\text{ is linearly independent.}
		\end{align*}
	\end{Lemma}
	\begin{proof}
		W.l.o.g. we assume that $\{k\in\{1,...,N\}:u_k\neq 0\} = \{1,...,l\}$, where $l\geq m$. We show linear independency straight forward with the result from the previous Lemma. Assume that  
		\begin{align*}
			0 = \sum\limits_{j = 0}^{r}\alpha_j\minimt{N}{t_j} = \sum\limits_{j = 0}^{r}\alpha_j\sum\limits_{k = 1}^{l}\frac{u_k}{1+t_j^2\lambda_k^2} \varphi_k =  \sum\limits_{k = 1}^{l}\sum\limits_{j = 0}^{r}\alpha_j\frac{u_k}{1+t_j^2\lambda_k^2} \varphi_k
		\end{align*}
		for some coefficients $\alpha_0,...,\alpha_r\in\R$. The orthonormal system $(\varphi_k)_{k= 1}^{l}$ is linearly independent, yielding
		\begin{align}\label{SecNormAppr:LinerarIndependencyEq}
			\llap{$\forall k\in\{1,...,l\}:\quad$}\sum\limits_{j = 	0}^{r}\alpha_j\frac{1}{1+t_j^2\lambda_k^2} = 0.
		\end{align}
		Let now $\mathbb{I}\subseteq \{1,...,l\}$, such that $\{\lambda_k^2\in\{\lambda_1^2,...,\lambda_N^2\}:k\in\mathbb{I}\} = \{\lambda_1^2,...,\lambda_l^2\}$ and $|\mathbb{I}| = m$. Then, condition \eqref{SecNormAppr:LinerarIndependencyEq} is equivalent to
		\begin{align*}
				\llap{$\forall k\in\mathbb{I}:\quad$}\sum\limits_{j = 	0}^{r}\alpha_j\frac{1}{1+t_j^2\lambda_k^2} = 0.
		\end{align*}
		Exploiting $r+1\leq m$ and the pairwise distinctness of both $(t_j)_{j=0}^r$ and $(\lambda_k^2)_{k\in\mathbb{I}}$, one deduces that $\alpha_0 = ... = \alpha_r = 0$. 
	\end{proof}
	The proof of Lemma \ref{Lm:LinInd} immediately reveals that the set $\{\minimt{N}{t_0},...,\minimt{N}{t_r}\}$ becomes linearly dependent as $r+1>m$. In practice, we observe two possible constellations. In the common case, where $u$ provides contributions from multiple basis vectors $\varphi_k$, solutions of \eqref{SecNormAppr:ShiftLapl} for different shift parameters $t_i$ indeed lead to an enrichment of $\Vr$, as long as $r$ is small enough. If the amount of non-zero Fourier-components of $u$ is rather small, we might be confronted with the case, where augmenting $r$ does not affect the dimension of $\Vr$ any further. In this case, however, enlargement of $\Vr$ is no longer necessary, as the following Theorem shows.
	\begin{Th}\label{Th:ExactRBNorm}
		Let $|\{\lambda_k^2\in\{\lambda_1^2,...,\lambda_N^2\}:u_k\neq 0\}| =: m\in\N$. If $r+1\geq m$, then the reduced basis interpolation norms \eqref{SecNormAppr:RBExtNorm}-\eqref{SecNormAppr:RBKNorm} coincide with the exact finite element interpolation norms from \eqref{SecInt:DiscIntNorms}, respectively.
	\end{Th}	
	\begin{proof}
		It suffices to validate the claim with respect to the Hilbert space interpolation norm for the case $r + 1 = m$. Let therefore $\mathbb{I}=\{i_0,...,i_{r}\}\subseteq\{1,...,N\}$, such that $\{\lambda_{i_0}^2,...,\lambda_{i_{r}}^2\} = \{\lambda_k^2\in\{\lambda_1^2,...,\lambda_N^2\}:u_k\neq 0\}$. Moreover, let
		\begin{align*}
			u = \bigoplus\limits_{i\in\mathbb{I}} u^{i}
		\end{align*}
		refer to the orthogonal decomposition of $u$ according to the corresponding  eigenspaces. Due to
		\begin{align*}
			\minimt{N}{t_j} = \sum\limits_{i\in\mathbb{I}} \frac{1}{1+t_j^2\lambda_i^2}u^i ,\rlap{\qquad$j = 0,...,r,$}
		\end{align*}
		and the regularity of the matrix $B\in\mathbb{R}^{(r+1)\times(r+1)}$ with $B_{kl} := (1+t_k^2\lambda_l^2)^{-1}$, we have that $(u^i)_{i\in\mathbb{I}}$ is a basis of $\Vr$, which is orthogonal by construction. We deduce $\lambda_{i_j}^2 = \mu_j^2$ for all $j = 0,...,r$. Direct computations reveal
		\begin{align*}
			\IntNormRB{u}{H}^2 = \sum\limits_{j=0}^{m-1}\mu_j^{2s}\scp{u}{\phi_j}{0}^2 = \sum\limits_{j=0}^{m-1}\lambda_{i_j}^{2s}\frac{\scp{u}{u^{i_j}}{0}^2}{\Norm{u^{i_j}}{0}^2}
			= \sum\limits_{j=0}^{m-1}\lambda_{i_j}^{2s}\Norm{u^{i_j}}{0}^2
			= \IntNormFEM{u}{H}^2.
		\end{align*}
	\end{proof}
	With exception of Theorem \ref{Th:ExactRBOp}, we only consider the case where $r+1\leq m$ for the rest of this paper, such that the dimension of $\Vr$ coincides with $r+1$. 
	\subsection{Computational aspects} The remainder of this section reviews the major ingredients to supply the reduced basis interpolation norms with a computationally applicable form. Before addressing this issue explicitly, we specify some further notation. Throughout what follows, for each $v\in\Vh$ we denote by $\CoefVec{v}\in\R^N$ its uniquely assigned coefficient vector, such that
	\begin{align*}
		v = \sum\limits_{k=1}^N (\CoefVec{v})_k b_k,
	\end{align*}	
	where $(b_k)_{k=1}^N$ denotes the finite element basis from \eqref{SecInt:DefMassStiffMat}. Moreover, we introduce the $s^{th}$-power of any symmetric matrix $Q\in\R^{l\times l}$, $l\in\N$, by diagonalization, i.e.,
	\begin{align*}
		Q^s := \Phi\Lambda^s\Phi^{-1},
	\end{align*}
	where $\Phi\in\R^{l\times l}$ denotes the matrix of eigenvectors of $Q$ and $\Lambda^s$ the involved diagonal matrix, containing the $s^{th}$-power of all corresponding eigenvalues. If $Q$ is also positive definite, we set
	\begin{align*}
		\llap{$\forall x,y\in\R^l:\quad$}\Norm{x}{Q}^2 :=x^TQx, \qquad \scp{x}{y}{Q} := x^TQy.
	\end{align*} 
		
	Based on these definitions, we follow the idea of POD \cite{KunischVolkwein} to propose an accurate procedure that computes \eqref{SecNormAppr:RBExtNorm} - \eqref{SecNormAppr:RBKNorm} for one and hence for all norms efficiently. The structure of $\|u\|_{E_r^s}$ and $\|u\|_{K_r^s}$ is less amenable to direct computations, since this would require quadrature rules on the unbounded domain $\R^+$. A more convenient setting can be provided by the eigensystem. Targeting at the computation of $\Vr$, we introduce the matrix
	\begin{align*}
		\widehat{V}_r := [\CoefVec{\minim{N}}(t_0),...,\CoefVec{\minim{N}}(t_r)] \in\R^{N\times(r+1)},
	\end{align*}
	whose $j^{th}$ column consists of the coefficient vector of $\minim{N}(t_j)$, i.e.,
	\begin{align}\label{SecNormAppr:LGS}
		\CoefVec{\minim{N}}(t_j) = (M + t_j^2A)^{-1}M\CoefVec{u}.
	\end{align}
	Thereupon, we introduce an orthonormal basis of $\Vr$ that is suitable for the study of our problem.
	\begin{Def}\label{Def:RBMatrix}
		The reduced basis matrix $V_r\in\R^{N\times(r+1)}$ is defined as the unique matrix that arises from Gram-Schmidt orthonormalization chronologically applied to the columns of $\widehat{V}_r$ with respect to the scalar product $\scp{\cdot}{\cdot}{M}$.
	\end{Def}
	\begin{Rem}\label{Rem:FirstColumnOfVr}
		The chronological performance of Gram-Schmidt orthonormalization in Definition \ref{Def:RBMatrix} yields that the first column of $V_r$ coincides with $\beta^{-1}\CoefVec{u}$, where $\beta := \Norm{u}{0}$.
	\end{Rem}
	The reduced basis matrix suggests a canonical basis on $\Vr$ by referring to the unique functions $b_{1}^r,...,b_{r}^r\in\Vr\subseteq\Vh$, whose assigned coefficient vectors coincide with the columns of $V_r$, i.e.,
	\begin{align*}
		V_r = \left[\CoefVec{b_{1}^r},...,\CoefVec{b_{r}^r}\right].
	\end{align*}
	Thereupon, we introduce for all $v_r\in\Vr$ its uniquely assigned coefficient vector $\CoefvecVr{v_r}\in\R^{r+1}$, such that
	\begin{align*}
		v_r = \sum\limits_{j=0}^r\left(\CoefvecVr{v_r}\right)_jb_{j}^r.
	\end{align*} 	
	There holds
	\begin{subequations}
		\begin{alignat}{2}
			&\scp{v_r}{w_r}{1} = \scp{\CoefVec{v_r}}{\CoefVec{w_r}}{A} = \CoefVec{v_r}^TA\CoefVec{w_r} = \CoefvecVr{v_r}^TV_r^TAV_r\CoefvecVr{w_r}, \label{SecNormAppr:Ar} \\			
			&\scp{v_r}{w_r}{0} = \CoefvecVr{v_r}^TV_r^TMV_r\CoefvecVr{w_r} = \CoefvecVr{v_r}^T\CoefvecVr{w_r}\label{SecNormAppr:Mr}
		\end{alignat}
	\end{subequations}
	for all $v_r, w_r\in\Vr$, where the last equality follows by the orthonormal property of $V_r$. Thereupon, we introduce the following definition.
	\begin{Def}\label{Def:ProcStiffMat}
		The projected stiffness matrix $A_r\in\R^{(r+1)\times(r+1)}$ is defined by
		\begin{align*}
			A_r := V_r^TAV_r.
		\end{align*}
	\end{Def}
	\begin{Th}\label{Th:CompOfRBNorm}
		Let $\CoefvecVr{e_1}\in\R^{r+1}$ denote the first unit vector and $\beta = \|u\|_{0}$. Then there holds
		\begin{align*}
			\|u\|_{H_{r}^s} = \beta\|\CoefvecVr{e_1}\|_{A_r^s}.
		\end{align*}
	\end{Th}
	The proof is postponed to Section \ref{SecOpAppr}. Theorem \ref{Th:CompOfRBNorm} highlights the beneficial structure of the proposed reduced basis algorithm. The arising problem size is of much smaller magnitude $r\ll N$, making direct computations of the eigensystem affordable.

	\section{Approximation of the operators}\label{SecOpAppr}
	All three reduced basis interpolation norms \eqref{SecNormAppr:RBExtNorm}-\eqref{SecNormAppr:RBKNorm} induce an operator, $\InducedRBOp{E}$, $\InducedRBOp{H}$, and $\InducedRBOp{K}$ on $\Vr$, which do all coincide up to $s$-dependent constants. In the same way as $\InducedRBOp{H}$ serves as reduced basis surrogate for the spectral finite element operator 
	\begin{align*}
		\InducedFEMOp{H} = \sum\limits_{k=1}^N\lambda_k^{2s}\scp{\cdot}{\varphi_k}{0}\varphi_k,
	\end{align*}
	$\InducedRBOp{E}$ can be interpreted as model order reduction of the generalized Dirichlet-to-Neumann map from \cite{HarmExtGeneralized}. The latter coincides with $\InducedOp{E}$, if $\Vnull$ and $\Veins$ are chosen appropriately. The purpose of $\InducedRBOp{K}$ is more of a technical than theoretical kind. We make use of the $\textrm{K}$-method as vital tool to proof convergence for one and hence for all three reduced basis operators. However, due to its computationally beneficial form, we stick to the equivalent spectral setting at first and refer to $\InducedRBOp{H}$ as our truth reduced basis approximation. In dependency of $u\in\Vh$, we state the essential definition of this paper.
	\begin{Def}[Reduced basis operator]\label{Def:RBOperator}
		For all $r\in\N$ we define the reduced basis operator $\InducedRBOp{H}$ of $\InducedFEMOp{H}$ as the induced operator of $\IntNormRB{\cdot}{H}$.
	\end{Def}
	\begin{Rem}
		In order to indicate its nonlinear nature, we write arguments of the reduced basis operator in brackets, i.e., $\InducedRBOp{H}(u)$ instead of $\InducedRBOp{H} u$.
	\end{Rem}
	As one can show, the matrix representation of $\InducedFEMOp{H}$, i.e., the unique matrix $\InducedFEMMatOp{H}\in\R^{N\times N}$, such that
	\begin{align*}
		\CoefVec{\InducedFEMOp{H}u} = \InducedFEMMatOp{H}\CoefVec{u}
	\end{align*}
	for all $u\in\Vh$, is given by $(M^{-1}A)^s$. It serves as matrix approximation of the original fractional operator $\InducedOp{H}$ in a finite element setting. Direct computations, however, are not feasible due to the typically large problem size. The proposed algorithm provides a remedy for this difficulty. In the following, we derive the $u$-dependent matrix representation $\InducedRBMatOp{H}$ of the reduced basis operator, providing the necessary information to carry out the actual computations.
	
	\begin{Th}\label{Th:IndSCPRB}
	For all $u\in \Vh$ there holds
		\begin{align}\label{SecOpAppr:NormRep}
			\IntNormRB{u}{H}^2 = \CoefVec{u}^TMV_rA_r^sV_r^TM\CoefVec{u}.
		\end{align}
	Moreover, the induced scalar product $\InducedRBScp{\cdot}{\cdot}{H}$ on $(\Vr,\IntNormRB{\cdot}{H})$ satisfies
	\begin{align}\label{SecOpAppr:ScpRep}
		\InducedRBScp{v_r}{w_r}{H} = \CoefVec{v_r}^TMV_rA_r^sV_r^TM\CoefVec{w_r}
	\end{align}
	for all $v_r,w_r\in\Vr$. The matrix representation of $\InducedRBOp{H}$ is given by
	\begin{align*}
		\InducedRBMatOp{H} = V_rA_r^sV_r^TM.
	\end{align*}
	\end{Th}
	\begin{proof}
		It suffices to prove \eqref{SecOpAppr:NormRep}. Recalling \eqref{SecNormAppr:RbEigenpairs}, we define
		\begin{align*}
			\Phi_r := \left[\CoefvecVr{\phi_0},...,\CoefvecVr{\phi_r}\right]\in\R^{(r+1)\times(r+1)}, \qquad \Lambda_r :=\textnormal{diag}(\mu_0^2,...,\mu_r^2) \in\R^{(r+1)\times(r+1)}.
		\end{align*}
		Equation \eqref{SecNormAppr:RbEigenpairs} combined with \eqref{SecNormAppr:Ar} and \eqref{SecNormAppr:Mr} yields
		\begin{align*}
			\Phi_r^T\Phi_r = I_r,\qquad \Phi_r^TA_r\Phi_r = \Lambda_r.
		\end{align*}
		Thus,
		\begin{align*}
			\CoefVec{u}^TMV_rA_r^sV_r^TM\CoefVec{u} = \CoefVec{u}^TM(V_r\Phi_r)\Lambda_r^s(V_r\Phi_r)^TM\CoefVec{u} = \sum_{j=0}^{r}\mu_j^{2s}\scp{u}{\phi_j}{0}^2 = \IntNormRB{u}{H}^2.
		\end{align*}
	\end{proof}
	We catch up on the postponed proof from the previous section.
	\begin{proof}[Proof of Theorem \ref{Th:CompOfRBNorm}]
		Follows directly from Theorem \ref{Th:IndSCPRB} and Remark \ref{Rem:FirstColumnOfVr} by utilizing the substitution $\CoefVec{u} = \beta V_r\CoefvecVr{e_1}$.
	\end{proof}
	
	$\InducedRBMatOp{H}$ serves as efficient approximation of $(M^{-1}A)^s$. Ahead of investigating its accuracy, we examine the arising computational costs, requiring knowledge of the map $r\mapsto (t_1,...,t_r)$ and its involved complexity. We address this problem adequately in Section \ref{Sec:ConvAna}, indicating for now that its essential computational expense amounts to finding a lower and upper bound for the spectrum of $M^{-1}A$. The overall complexity has to be regarded from two different perspectives, the so-called \textit{offline} and \textit{online} phase, and depends on the particular problem.

	At first, we consider evaluations of type $u\mapsto \InducedRBOp{H}(u)$, which are of nonlinear character. The underlying offline phase encompasses computations of the spectral bounds as a one-time investment. The online phase has to be performed for each argument separately. It incorporates computations of $r$ finite element solutions $\minimt{N}{t_j}$ of the original, expensive problem size, followed by orthonormalization of $\widehat{V}_r\in\R^{N\times (r+1)}$ to obtain the reduced basis matrix $V_r$. Furthermore, the projected matrix $A_r = V_r^TAV_r$ has to be established in order to determine its eigensystem. The assembly of first $A_r^s$ and second of $V_rA_r^sV_r^TM\CoefVec{u}$ completes the computations. Despite its nonlinear nature, the savings gained substantially outweigh the arising inaccuracy, if $r$ is of moderate size.

	The offline-online decomposition is of particular interest, if we target at approximations of type $s\mapsto \InducedRBOp{H}(u)$ for fixed $u$ and several values of $s\in(0,1)$. In this case, the online phase breaks down to the assembly of $A_r^s$ and $V_rA_r^sV_r^TM\CoefVec{u}$, while all remaining computations are the matter of a one-time investment within the offline phase.
	
	Several properties of the reduced basis interpolation norms also apply to Definition \ref{Def:RBOperator}. Exemplary, the operator counterpart of Theorem \ref{Th:ExactRBNorm} is procured in the following.
	\begin{Th}\label{Th:ExactRBOp}
		Let $|\{\lambda_k^2\in\{\lambda_1^2,...,\lambda_N^2\}:u_k\neq 0\}| =: m\in\N$. If $r+1\geq m$, then there holds
		\begin{align*}
			\InducedRBOp{H}(u) = \InducedFEMOp{H}u.
		\end{align*}
	\end{Th}
	\begin{proof}
		Follows the very same arguments as the proof of Theorem \ref{Th:ExactRBNorm}.
	\end{proof}

	\section{Convergence Analysis}\label{Sec:ConvAna}
	The goal of this section is to specify the choice of snapshots in Definition \ref{Def:Method} to gain optimal convergence properties. We affirm that there exists a tuple of positive numbers $t_1,...,t_r$, naturally arising from the analysis, such that exponential decay in the error for both norm and operator action is obtained. While computations are carried out in the spectral setting, the approach involving the $\textrm{K}$-functional turns out to provide a more beneficial environment for the analysis. We adopt Theorem \ref{Th:IndSCPRB} in this context and take advantage of the, up to a multiplicative constant, interchangeable role of the reduced basis scalar products. Before going further into detail, we investigate two fundamental definitions that are based on \cite{ZolotarevCollectedWorks} and \cite{ZolotarevProbGonchar}, see also \cite{Oseledets}, involving the theory of elliptic integrals and Jacobi elliptic functions, see \cite[Section 16 \& 17]{HandbookOfMathFunc}.	
	\begin{Def}\label{Def:ZolotarevPoints}
		Let $\delta\in(0,1)$. For each $r\in\N$ we define the Zolotar\"ev points $\mathcal{Z}_1,...,\mathcal{Z}_r$ on the interval $[\delta,1]$ by
		\begin{align*}
			\mathcal{Z}_j := \DeltaAmpl\left(\frac{2(r-j)+1}{2r}\mathcal{K}(\delta'),\delta'\right), \rlap{\qquad$j = 1,...,r,$}
		\end{align*}
		where $\DeltaAmpl(\theta,k)$ denotes the Jacobi elliptic function, $\mathcal{K}(k)$ the corresponding elliptic integral of first kind with elliptic modulus $k$, and $\delta':= \sqrt{1-\delta^2}$.
	\end{Def}
	\begin{Def}\label{Def:TransformedZolotarevPoints}
		Let $0<a<b\in\R^+$. For each $r\in\N$ we define the transformed Zolotar\"ev points $\widehat{\mathcal{Z}}_1,...,\widehat{\mathcal{Z}}_r$ on $[a,b]$ by
		\begin{align*}
			\widehat{\mathcal{Z}}_j := b\mathcal{Z}_j, \rlap{\qquad$j = 1,...,r$,}
		\end{align*}
		where $\mathcal{Z}_1,...,\mathcal{Z}_r$ refer to the Zolotar\"ev points on $\left[\frac{a}{b},1\right]$.
	\end{Def}
	As shown in the further course of action, the transformed Zolotar\"ev points on $[\lambda_N^{-2},\lambda_1^{-2}]$ turn out to be perfectly tailored for our reduced basis strategy. We recap that $\lambda_1^2$ and $\lambda_N^2$ refer to the minimal and maximal eigenvalue of of the discrete operator and agree on the following nomenclature.
	\begin{Def}\label{Def:ZolotarevSpace}
		A reduced space $\Vr = \Span\{\minimt{N}{t_0},...,\minimt{N}{t_r}\}\subseteq \Vh$ is called Zolotar\"ev space, if and only if
		there exist two constants $\lambda_L^{2},\lambda_U^2\in\R^+$ with
		\begin{align*}
			\lambda_L^2\leq \lambda_1^2,	\qquad\qquad \lambda_U^2 \geq \lambda_N^2,	
		\end{align*} 
			such that the squared snapshots $t_1^2,...,t_r^2$ coincide with the transformed Zolotar\"ev points on $\sigma^{inv} := [\lambda_U^{-2}, \lambda_L^{-2}]$. We call $\sigma := [\lambda_L^2,\lambda_U^2]$ the spectral interval of $\Vr$.
	\end{Def}
	
	\subsection{Error of the reduced basis interpolation norms}\label{SecConvAna:SubSecIntNormError}
	We specify some further notation. By $a\preceq b$ we mean that there exists a constant $C\in\R^+$, independent of $a$, $b$, and $r$, such that $a \leq  Cb$. Along with this premiss, we prove that Definition \ref{Def:ZolotarevSpace} provides an optimal choice for the reduced space. 
	\begin{Th}[Exponential convergence of the reduced basis interpolation norms]\label{Th:NormError}
		Let $u\in \Vh$ and $\Vr\subseteq \Vh$ a Zolotar\"ev space with $\sigma = [\lambda_L^2,\lambda_U^2]$ and $\delta = \nicefrac{\lambda_L^2}{\lambda_U^2}$. Then there holds
		\begin{align}\label{SecConvAna:NormError}
			\llap{$\exists\hspace{0.05cm} C\in\R^+:$\quad}	0\leq \IntNormRB{u}{K}^2 - \IntNormFEM{u}{K}^2  \preceq e^{-Cr}\Norm{u}{1}^2.
		\end{align}
		The constant $C$ only  depends on $\delta$ and satisfies
		\begin{align*}
			C(\delta) = \mathcal{O}\left(\frac{1}{\ln\left(\frac{1}{\delta}\right)}\right),\quad\text{as }\delta\to 0.
		\end{align*}
		Its precise value coincides with $2C^*$, where $C^*$ refers to the constant from Remark \ref{Rem:AsympBehavC}.
	\end{Th} 
	Verifying the claim of Theorem \ref{Th:NormError} is challenging, which is why we conduct the proof in several steps. Under the prescribed assumptions, we start with the first inequality of \eqref{SecConvAna:NormError}.
	\begin{Lemma}
		There holds
		\begin{align*}
			\IntNormRB{u}{K} \geq \IntNormFEM{u}{K}.
		\end{align*}
	\end{Lemma}
	\begin{proof}
		The relation $\Vr\subseteq \Vh$ immediately reveals
		\begin{align*}
			\llap{$\forall t\in\R^+:$\quad} K_r^2(t;u) &= \inf_{v_r\in\Vr} \|u- v_r\|_{0}^2 + t^2\|v_r\|_{1}^2\\
		&\geq \inf_{v\in \Vh} \|u-{v}\|_{0}^2 + t^2\|{v}\|_{1}^2 = K^2(t;u).
		\end{align*}
		Hence,
		\begin{align*}
			\IntNormRB{u}{K}^2 = \int_{0}^{\infty}t^{-2s-1}K_r^2(t;u)\,dt \geq \int_{0}^{\infty}t^{-2s-1}K^2(t;u)\,dt = \IntNormFEM{u}{K}^2.
		\end{align*}
	\end{proof}
	Due to 
	\begin{align}\label{SecConvAna:NormErrorRepresentation}
		\IntNormRB{u}{K}^2 - \IntNormFEM{u}{K}^2 = \int_{0}^{\infty}t^{-2s-1}\left(K_r^2(t;u) - K^2(t;u)\right)\,dt,
	\end{align}
	the error of the norms can be traced back to the error of the $\textrm{K}$-functionals. To make matters precise, we conduct some technical preparations.
	\begin{Def}
		For all $t\in\R^+$ we denote by $\minimt{r}{t}\in\Vr$ the unique minimizer of $K_r^2(t;u)$.
	\end{Def}
	\begin{Rem}\label{Rem:ShiftedLaplaceRb}
		Similarly to $\minimt{N}{t}$, utilizing Euler-Lagrange formalism, the minimizer $\minimt{r}{t}\in\Vr$ is the unique solution of the variational problem
		\begin{align*}
			\llap{$\forall w_r\in\Vr:$\quad} \scp{\minimt{r}{t}}{w_r}{0} + t^2\scp{\minimt{r}{t}}{w_r}{1} = \scp{u}{w_r}{0},
		\end{align*}
		or equivalently, $\CoefvecVr{\minim{r}}(t)\in\R^{r+1}$ solves the linear system of equations
		\begin{align*}
			(I_r + t^2A_r)\CoefvecVr{\minim{r}}(t) = V_r^TM\CoefVec{u},
		\end{align*}
		where $I_r\in\R^{(r+1)\times(r+1)}$ represents the identity matrix.
	\end{Rem}
	
	\begin{Lemma}\label{Lm:TechnicalDiffK}
		For all $t\in\R^+$ and $w_r\in\Vr$ there holds
		\begin{align*}
			\|u-w_r\|_{0}^2 + t^2\|w_r\|_{1}^2 - \|u-\minimt{N}{t}\|_{0}^2 - t^2&\|\minimt{N}{t}\|_{1}^2 =  \\ &\|\minimt{N}{t} - w_r\|_{0}^2 + t^2\|\minimt{N}{t} - w_r\|_{1}^2.
		\end{align*}
	\end{Lemma}
	\begin{proof}
		For convenience, we omit the dependency in $t$ in consecutive elaborations. One observes
		\begin{align*}
			\|u-w_r\|_{0}^2 - \|u-\minim{N}\|_{0}^2 &= \scp{u-w_r}{u-w_r}{0} - \scp{u-\minim{N}}{u-\minim{N}}{0} \\
			&= -2\scp{u}{w_r}{0} + \|w_r\|_{0}^2 + 2\scp{u}{\minim{N}}{0} - \|{\minim{N}}\|_{0}^2.
		\end{align*}
		We define the bilinear form $a(u,w) := \scp{u}{w}{0} + t^2 \scp{u}{w}{1}$ on $\Vh\times\Vh$. Due to \eqref{SecNormAppr:ShiftLapl}, there holds
		\begin{align*}
			\|u-w_r\|_{0}^2 - \|u-\minim{N}\|_{0}^2 = -2a(\minim{N},w_r) + \|w_r\|_{0}^2 + 2a(\minim{N},\minim{N}) - \|\minim{N}\|_{0}^2.
		\end{align*}
		Thus,
		\begin{align*}
			\|u-w_r\|_{0}^2 + t^2\|w_r\|_{1}^2 - \|u-\minim{N}\|_{0}^2 - t^2\|\minim{N}\|_{1}^2
			&= a(w_r,w_r) - 2a(\minim{N},w_r) + a(\minim{N},\minim{N})\\
			&= a(w_r -\minim{N},w_r -\minim{N}) \\
			&=\|\minim{N} - w_r\|_{0}^2 + t^2\|\minim{N} - w_r\|_{1}^2.
		\end{align*}
	\end{proof}
	The accuracy of $K_r^2(t;u)$ rests upon the approximation quality of the minimizer $\minimt{r}{t}\approx\minimt{N}{t}$, as the following result shows.
	\begin{Cor}\label{Cor:DiffKDiffInf}
		For all $t\in\R^+$ there holds
		\begin{align}\label{SecConcAna:DiffK=DiffInf}
			K_r^2(t;u) - K^2(t;u) = \| \minimt{N}{t} - \minimt{r}{t}\|_{0}^2 + t^2\|\minimt{N}{t} - \minimt{r}{t}\|_{1}^2.
		\end{align}
	\end{Cor}	
	\begin{proof}
		Follows directly from Lemma \ref{Lm:TechnicalDiffK} with $w_r = \minimt{r}{t}$.
	\end{proof}
	Dealing with \eqref{SecConcAna:DiffK=DiffInf} is challenging. We derive an upper bound for the error which turns out to be more amenable to analytical considerations.
	\begin{Cor}\label{Cor:UpperBoundDiffK}
		For all $t\in\R^+$ and $w_r\in\Vr$ there holds
		\begin{align}\label{SecConvAnaUpperBoundKFunc}
			K_r^2(t;u) - K^2(t;u) \leq \|\minimt{N}{t} - w_r\|_{0}^2 + t^2\|\minimt{N}{t} - w_r\|_{1}^2.
		\end{align}
	\end{Cor}
	\begin{proof}
		Due to the minimization property of $\minimt{r}{t}$, there holds
		\begin{align*}
			K_r^2(t;u) - K^2(t;u) \leq \|u-w_r\|_{0}^2 + t^2\|w_r\|_{1}^2 - \|u-\minimt{N}{t}\|_{0}^2 - t^2\|\minimt{N}{t}\|_{1}^2
		\end{align*}
		for all  $t\in\R^+$ and $w_r\in\Vr$. Utilizing Lemma \ref{Lm:TechnicalDiffK} concludes the proof.
	\end{proof}
	In the subsequent, we aim to choose
	\begin{align}\label{SecConvAna:wr}
		w_r = \sum\limits_{j=0}^{r}\alpha_j\minimt{N}{t_j}\in\Vr
	\end{align}
	from Corollary \ref{Cor:UpperBoundDiffK} in a clever way, such that the upper bound in \eqref{SecConvAnaUpperBoundKFunc} becomes small. The idea of how to choose its coefficients $\alpha_j$ emerges from the following investigation. 
	\begin{Th}\label{Thm:KFuncInterpBound}
		For all $\alpha_0,...,\alpha_r\in\R$ there holds
		\begin{align*}
			K_r^2(t;u) - K^2(t;u) \leq \sum\limits_{k=1}^N (1 + t^2\lambda_k^2)\Bigg(\frac{1}{1 + t^2\lambda_k^2} - \sum\limits_{j = 0}^r \alpha_j\frac{1}{1 + t_j^2\lambda_k^2}\Bigg)^2u_k^2.
		\end{align*}
	\end{Th}
	\begin{proof}
		Due to Corollary \ref{Cor:UpperBoundDiffK}, there holds for any function of type \eqref{SecConvAna:wr}
		\begin{align*}
			K_r^2(t;u) - K^2(t;u) \leq \|\minimt{N}{t} - w_r\|_{0}^2 + t^2\|\minimt{N}{t} - w_r\|_{1}^2.
		\end{align*}
		Spectral decomposition combined with Lemma \ref{Lm:InfRepr} yields	
		\begin{align*}
			\minimt{N}{t} - w_r &= \sum\limits_{k=1}^N\frac{u_k}{1 + t^2\lambda_k^2}\varphi_k - \sum\limits_{j = 0}^r \alpha_j\sum\limits_{k=1}^N\frac{u_k}{1 + t_j^2\lambda_k^2}\varphi_k \\
			&= \sum\limits_{k=1}^N\Bigg(\frac{1}{1 + t^2\lambda_k^2} - \sum\limits_{j = 0}^r \alpha_j\frac{1}{1 + t_j^2\lambda_k^2}\Bigg)u_k\varphi_k.
		\end{align*}	
		Based on the orthogonal property of $(\varphi_k)_{k=1}^N$, we observe
		\begin{align*}	
			\|\minimt{N}{t} - w_r\|_{1}^2 &= \sum\limits_{k=1}^N\Big\|\Bigg(\frac{1}{1 + t^2\lambda_k^2} - \sum\limits_{j = 0}^r \alpha_j\frac{1}{1 + t_j^2\lambda_k^2}\Bigg)u_k\varphi_k\Big\|_{1}^2 \\
			&= \sum\limits_{k=1}^N \lambda_k^2\Bigg(\frac{1}{1 + t^2\lambda_k^2} - \sum\limits_{j = 0}^r \alpha_j\frac{1}{1 + t_j^2\lambda_k^2}\Bigg)^2u_k^2.
		\end{align*}
		Computations in the $0$-norm can be concluded analogously, proving the claim.
	\end{proof}
	Theorem \ref{Thm:KFuncInterpBound} reveals that \eqref{SecConvAna:wr} has to be chosen in a way, such that for all $\lambda_1,...,\lambda_N$, or more generally, for all values of $\lambda\in[\lambda_1,\lambda_N]$, the difference
	\begin{align}\label{SecConvAna:RationalError}
		\left(\frac{1}{1 + t^2\lambda^2} - \sum\limits_{j = 0}^r {\alpha}_j\frac{1}{1 + t_j^2\lambda^2}\right)
	\end{align}
	becomes small. Typically, neither $\lambda_1$ nor $\lambda_N$ are known a-priori, which is why we consider \eqref{SecConvAna:RationalError} with respect to the spectral interval from Definition \ref{Def:ZolotarevSpace}, admitting $\lambda\in[\lambda_L,\lambda_U]$ instead. Any possible bound of \eqref{SecConvAna:RationalError} then trivially also holds on $[\lambda_1,\lambda_N]$.
	
	In the further course of action, we derive two different candidates for the coefficients $(\alpha_j)_{j=0}^r$ in dependency of $t$. The first one ensures that \eqref{SecConvAna:RationalError} becomes small for $t\geq1$, while the second achieves the same as $t<1$. To this extent, we make a first ansatz and set $\alpha_0 = 0$. The latter coefficients are determined by means of a rational interpolation problem, which we inquire in the subsequent Lemma.
	\begin{Lemma}\label{Lm:RatInt}
		Assume that $\kappa\in\R^+, \kappa_1,...,\kappa_r\in \sigma^{inv}$, and $\kappa_i\neq\kappa_j$ for $i\neq j$. Consider the space $\mathcal{R}$ of all rational functions $R$, which admit a representation
		\begin{align*}
			R(x) = \sum\limits_{j=1}^r\alpha_j\frac{1}{1+\kappa_jx}
		\end{align*}
		for coefficients $\alpha_1,...,\alpha_r\in\R$. Further define
		\begin{align*}
			g_{\kappa}(x) := \frac{1}{1+\kappa x}.
		\end{align*}
		Then the solution of the rational interpolation problem: Find $q\in\mathcal{R}$, such that
		\begin{align}\label{LmProof:RatInt}
			\llap{$\forall j\in\{1,...,r\}:\quad$} q\left(\frac{1}{\kappa_j}\right) = g_{\kappa}\left(\frac{1}{\kappa_j}\right)
		\end{align}
		satisfies
		\begin{align}\label{LmProof:ErrorProduct}
			\llap{$\forall x\in\sigma:$\quad}|g_{\kappa}(x)-q(x)| \leq \frac{1}{1+\kappa x}\prod\limits_{j=1}^{r}\left|\frac{1 - \kappa_jx}{1 + \kappa_jx}\right|.
		\end{align}
	\end{Lemma}
	\begin{proof}
		Let
		\begin{align*}
			q(x) = \sum\limits_{j = 1}^r{\hat{\alpha}}_j\frac{1}{1+\kappa_jx}
		\end{align*}
		denote the unique solution of \eqref{LmProof:RatInt}. Then there holds
		\begin{align}\label{IntError}
			g_{\kappa}(x) - q(x) = \frac{1}{1+\kappa x} - \sum\limits_{j = 1}^{r}{\hat{\alpha}}_j\frac{1}{1+\kappa_j x} = \frac{p(x)}{(1+\kappa x)\prod\limits_{j = 1}^{r}(1 + \kappa_j x)}
		\end{align}
		for a suitable polynomial $p$ of degree $r$. The interpolation property yields
		\begin{align*}
			\llap{$\forall j\in\{1,...,r\}:\quad$} p\left(\frac{1}{\kappa_j}\right) = 0 .
		\end{align*}
		The fundamental Theorem of algebra affirms the existence of a constant $c \in\R$, such that
		\begin{align*}
			p(x) = c\prod\limits_{j=1}^{r}(1 - \kappa_jx).
		\end{align*}
		The constant $c$ can be further specified. Multiplying \eqref{IntError} by $(1+\kappa x)$ and setting $x = -\frac{1}{\kappa}$ yields
		\begin{align*}
			1 = c\prod\limits_{j=1}^r\frac{1+\frac{\kappa_j}{\kappa}}{1-\frac{\kappa_j}{\kappa}}, \quad \text{and hence,} \quad c = \prod\limits_{j=1}^{r}\frac{(\kappa - \kappa_j)}{(\kappa + \kappa_j)}.
		\end{align*}
		All together, we obtain for all $x\in\sigma$
		\begin{align*}
			|g_{\kappa}(x) - q(x)| = \frac{1}{1+\kappa x}\prod\limits_{j=1}^{r}\left|\frac{(\kappa - \kappa_j)}{(\kappa+\kappa_j)}\frac{(1 - \kappa_j x)}{(1 + \kappa_j x)}\right| \leq \frac{1}{1+\kappa x}\prod\limits_{j=1}^{r}\left|\frac{1 - \kappa_j x}{1 + \kappa_j x}\right|.
		\end{align*}
	\end{proof}
	Minimizing the maximal deviation of the upper bound in \eqref{LmProof:ErrorProduct} leads to a min-max problem of the following kind: Find $\kappa_1,...,\kappa_r\in \sigma^{inv}$, such that
	\begin{align}\label{SecConvAna:MinMax}
		\min\limits_{\theta_1,...,\theta_r\in\sigma^{inv}}\max\limits_{x\in\sigma} \prod\limits_{j=1}^{r}\left|\frac{1 - \theta_j x}{1 + \theta_j x}\right| = \max\limits_{x\in \sigma}\prod\limits_{j=1}^{r}\left|\frac{1 - \kappa_j x}{1 + \kappa_j x}\right|.
	\end{align}
	Closely related problems have been investigated in \cite{ZolotarevCollectedWorks} and \cite{ZolotarevProbGonchar}. We summarize the essential results. Consider the slightly modified problem: Find $\kappa_1,...,\kappa_r\in [\delta,1]$, $\delta = \nicefrac{\lambda_L^2}{\lambda_U^2}$, such that
	\begin{align}\label{SecConvAna:ZolotarevMinMax}
		\min\limits_{\theta_1,...,\theta_r\in[\delta,1]}\max\limits_{x\in[\delta,1]} \prod\limits_{j=1}^{r}\left|\frac{x-\theta_j}{x + \theta_j}\right| = \max\limits_{x\in[\delta,1]}\prod\limits_{j=1}^{r}\left|\frac{x - \kappa_j}{x+\kappa_j}\right|.
	\end{align}
	Zolotar\"ev and Gonchar showed that its unique solution is given by the Zolotar\"ev points $\mathcal{Z}_1,...,\mathcal{Z}_r$ on $[\delta,1]$. They further approved that there exists a positive constant $C$, depending on $\delta$ only, such that
	\begin{align}\label{SecConvAna:ZoloProd}
		\max\limits_{x\in[\delta,1]}\prod\limits_{j=1}^{r}\left|\frac{x - \mathcal{Z}_j}{x+\mathcal{Z}_j}\right|\preceq e^{-Cr}.
	\end{align}
	The product in \eqref{SecConvAna:ZoloProd}, considered as function in $x$, features $r+1$ points of alternance and has the least deviation from zero on $[\delta,1]$ among all functions of this form. We set results from problem \eqref{SecConvAna:ZolotarevMinMax} in correspondence with \eqref{SecConvAna:MinMax}.
	\begin{Th}\label{Th:OriginalMinMaxProblem}
		The unique solution $\kappa_1,...,\kappa_r$ of problem \eqref{SecConvAna:MinMax} satisfies $\kappa_j = \widehat{\mathcal{Z}}_j$ for all $j=1,...,r$, where $\widehat{\mathcal{Z}}_j$ denotes the $j^{th}$ transformed Zolotar\"ev point on $\sigma^{inv}$.
	\end{Th}
	\begin{proof}
		Consider the linear transformation $\Psi:[\delta,1]\longrightarrow\sigma^{inv}$ with $\Psi(x) := \lambda_L^{-2}x$. Direct computations, based on the results of \cite{ZolotarevCollectedWorks}, reveal 
		\begin{align*}
			\min\limits_{\theta_1,...,\theta_r\in\sigma^{inv}}\max\limits_{x\in\sigma} \prod\limits_{j=1}^{r}\left|\frac{1-\theta_jx}{1 + \theta_j x}\right| &= \min\limits_{\theta_1,...,\theta_r\in\sigma^{inv}}\max\limits_{x\in\sigma^{inv}} \prod\limits_{j=1}^{r}\left|\frac{x-\theta_j}{x + \theta_j}\right| \\
			&=\min\limits_{\theta_1,...,\theta_r\in[\delta,1]}\max\limits_{x\in[\delta,1]} \prod\limits_{j=1}^{r}\left|\frac{\Psi(x)-\Psi(\theta_j)}{\Psi(x) + \Psi(\theta_j)}\right| \\
			&=\max\limits_{x\in[\delta,1]} \prod\limits_{j=1}^{r}\left|\frac{x-\mathcal{Z}_j}{x + \mathcal{Z}_j}\right| \\
			&= \max\limits_{x\in[\delta,1]} \prod\limits_{j=1}^{r}\left|\frac{\Psi(x)-\Psi(\mathcal{Z}_j)}{\Psi(x) + \Psi(\mathcal{Z}_j)}\right| \\
			&= \max\limits_{x\in\sigma^{inv}} \prod\limits_{j=1}^{r}\left|\frac{x-\widehat{\mathcal{Z}}_j}{x + \widehat{\mathcal{Z}}_j}\right| = \max\limits_{x\in\sigma}\prod\limits_{j=1}^{r}\left|\frac{1-\widehat{\mathcal{Z}}_jx}{1 + \widehat{\mathcal{Z}}_j x}\right|.
		\end{align*}
	\end{proof}
	Within the proof of Theorem \ref{Th:OriginalMinMaxProblem}, we have additionally shown that 
	\begin{align*}
		\max\limits_{x\in\sigma} \prod\limits_{j=1}^{r}\left|\frac{1-\widehat{\mathcal{Z}}_jx}{1 + \widehat{\mathcal{Z}}_j x}\right| = \max\limits_{x\in[\delta,1]} \prod\limits_{j=1}^{r}\left|\frac{x-\mathcal{Z}_j}{x + \mathcal{Z}_j}\right|,
	\end{align*}
	which, due to \eqref{SecConvAna:ZoloProd}, immediately reveals the following result.
	\begin{Cor}\label{Cor:ZolotProdConv}
		Let $\widehat{\mathcal{Z}}_1,...,\widehat{\mathcal{Z}}_r$ denote the transformed Zolotar\"ev points on $\sigma^{inv}$. Then there holds
		\begin{align}\label{SecConvAna:MaxDev}
			\llap{$\exists\hspace{0.05cm} C\in\R^+:$\quad}\max\limits_{x\in\sigma} \prod\limits_{j=1}^{r}\left|\frac{1-\widehat{\mathcal{Z}}_jx}{1 + \widehat{\mathcal{Z}}_j x}\right|\preceq e^{-Cr}.
		\end{align}
	\end{Cor}
	\begin{Rem}\label{Rem:AsympBehavC}
		As shown in \cite{MedovikovAndLebedev}, constants from Corollary \ref{Cor:ZolotProdConv} can be further specified. More precisely, the maximal deviation in \eqref{SecConvAna:MaxDev} can be bounded by means of
		\begin{align*}
			\max\limits_{x\in\sigma} \prod\limits_{j=1}^{r}\left|\frac{1-\widehat{\mathcal{Z}}_jx}{1 + \widehat{\mathcal{Z}}_j x}\right| \leq 2e^{-C^*r},
		\end{align*}
		with
		\begin{align*}
			C^* := \frac{\pi \mathcal{K}(\mu_1)}{4\mathcal{K}(\mu)},\qquad\quad \mu := \left(\frac{1-\sqrt{\delta}}{1+\sqrt{\delta}}\right)^2,\qquad\quad \mu_1 := \sqrt{1-\mu^2},
		\end{align*}
		and $\mathcal{K}$ the elliptic integral from Definition \ref{Def:ZolotarevPoints}. The following asymptotic formulas are known to hold,
		\begin{align*}
			\mathcal{K}(\mu) \approx \frac{1}{2}\ln\left(\frac{16}{1-\mu}\right), \qquad \mathcal{K}(\mu_1)\approx \frac{\pi}{2}, \quad\text{as }\mu\to 1,
		\end{align*} 
		see e.g., \cite[Section 17]{HandbookOfMathFunc}. This yields the asymptotic behaviour of $C^*$ in dependency of $\delta$,
		\begin{align*}
			C^*(\delta) \approx \frac{1}{\ln\left(\frac{1}{\delta}\right)}, \quad\text{as }\delta\to0.
		\end{align*} 
		Along with the choice $\lambda_L := \lambda_1$ and $\lambda_U := \lambda_N$, this reveals that the constant $C^*$ only deteriorates at logarithmical rate as the condition number $\nicefrac{\lambda_N^2}{\lambda_1^2} = \delta^{-1}$ increases. The number of solves required to achieve a prescribed precision $\varepsilon > 0$ behaves like
		\begin{align*}
			r = \mathcal{O}(\ln(\varepsilon)\ln(\delta)).
		\end{align*} 
	\end{Rem}
	We eventually return to the original problem of interest in \eqref{SecConvAna:RationalError}.
	\begin{Lemma}\label{Lm:IntBound}
		Denote by $\hat{\alpha}_1,...,\hat{\alpha}_r\in\R$ the coefficients of the unique solution $q\in\mathcal{R}$ from the rational interpolation problem \eqref{LmProof:RatInt} with $\kappa := t^2$ for some $t\in\R^+$ and $\kappa_j := t_j^2$. Let $C^*$ denote the constant from Remark \ref{Rem:AsympBehavC}. Then there holds
		\begin{align}\label{SecConvAna:DiskDiffRatInt}
			\left|\frac{1}{1 + t^2\lambda_k^2} - \sum\limits_{j = 1}^r \hat{\alpha}_j\frac{1}{1 + t_j^2\lambda_k^2}\right|\preceq \frac{e^{-C^*r}}{1+t^2\lambda_k^2}, \rlap{\qquad$k=1,...,N.$}
		\end{align}
	\end{Lemma}
	\begin{proof}
		According to Lemma \ref{Lm:RatInt}, Corollary \ref{Cor:ZolotProdConv}, and Remark \ref{Rem:AsympBehavC}, there holds
		\begin{align*}
			\left|\frac{1}{1 + t^2\lambda^2} - \sum\limits_{j = 1}^r \hat{\alpha}_j\frac{1}{1 + t_j^2\lambda^2}\right| \leq \frac{1}{1+t^2 \lambda^2}\prod\limits_{j=1}^{r}\left|\frac{1 - t_j^2 \lambda^2}{1 + t_j^2 \lambda^2}\right| \preceq \frac{e^{-C^*r}}{1+t^2\lambda^2}
		\end{align*}
		for all $\lambda\in[\lambda_L,\lambda_U]$. By construction, there holds $\lambda_k\in[\lambda_L,\lambda_N]$ for all $k = 1,...,N$, which is why \eqref{SecConvAna:DiskDiffRatInt} is valid.
	\end{proof}
	Assembling results from above finally enables us to derive an upper bound for the error in the reduced basis $\textrm{K}$-functional.
	\begin{Th}\label{Th:ErrorKFunc}
		Let $C^*$ denote the constant from Remark \ref{Rem:AsympBehavC}. Then there holds for all $t\in\R^+$
		\begin{align*}
			K_r^2(t;u) - K^2(t;u) \preceq e^{-2C^*r}\sum\limits_{k = 1}^N\frac{1}{1+t^2\lambda_k^2}u_k^2.
		\end{align*}
	\end{Th}
	\begin{proof}
		Theorem \ref{Thm:KFuncInterpBound} combined with Lemma \ref{Lm:IntBound} yields
		\begin{align*}
			K_r^2(t;u) - K^2(t;u) &\preceq \sum\limits_{k = 1}^N(1+t^2\lambda_k^2)\left(\frac{e^{-C^*r}}{1+t^2\lambda_k^2}\right)^2u_k^2 = e^{-2C^*r}\sum\limits_{k=1}^N\frac{1}{1+t^2\lambda_k^2}u_k^2.
		\end{align*}
	\end{proof}
	As it turns out in the further course of action, the upper bound derived in Theorem \ref{Th:ErrorKFunc} is only sharp enough in case of $t\geq1$, but not as $t<1$. We overcome this inconvenience by subtle adjustments of the interpolation problem \eqref{LmProof:RatInt}, such that its arising solution leads to the desired properties. 
	\begin{Th}\label{Th:ErrorKFunc2}
		Let $C^*$ denote the constant from Remark \ref{Rem:AsympBehavC}. Then there holds for all $t\in\R^+$
		\begin{align*}
			K_r^2(t;u) - K^2(t;u) \preceq e^{-2C^*r}\sum\limits_{k = 1}^N\frac{t^4\lambda_k^4}{1+t^2\lambda_k^2}u_k^2.
		\end{align*}
	\end{Th}
	\begin{proof}
		In analogy to Lemma \ref{Lm:RatInt}, we consider the following rational interpolation problem: For $\kappa\in\R^+$ and $\kappa_1,...,\kappa_r\in\sigma^{inv}$ pairwise distinct, find $q\in\hat{\mathcal{R}}$, such that
		\begin{align*}
			\llap{$\forall j\in\{1,...,r\}:$\quad} q\left(\frac{1}{\kappa_j}\right) &= g_{\kappa}\left(\frac{1}{\kappa_j}\right), \\
			q(0) &= g_{\kappa}(0),
		\end{align*}
		where $\hat{\mathcal{R}}$ denotes the linear span of $\mathcal{R}$ enriched by constant functions. We remark that $\alpha_0 = 0$ is no longer constrained. Similarly to the proof of Lemma \ref{Lm:RatInt}, one affirms that
		\begin{align*}
			|g_{\kappa}(x) - q(x)| = \frac{\kappa x}{1+\kappa x}\prod\limits_{j=1}^{r}\left|\frac{(\kappa - \kappa_j)}{(\kappa+\kappa_j)}\frac{(1 - \kappa_j x)}{(1 + \kappa_j x)}\right| \leq \frac{\kappa x}{1+\kappa x}\prod\limits_{j=1}^{r}\left|\frac{1 - \kappa_j x}{1 + \kappa_j x}\right|.
		\end{align*}
		The transformed Zolotar\"ev points on $\sigma^{inv}$ ensure exponential convergence of the product. Proceeding in the same manner as before, one concludes the proof.
	\end{proof}
	
	We have now all the required tools to prove exponential convergence of the reduced basis interpolation norms in the $\textrm{K}$-setting.
	\begin{proof}[Proof of Theorem \ref{Th:NormError}]
		Theorem \ref{Th:ErrorKFunc} and \ref{Th:ErrorKFunc2} together with \eqref{SecConvAna:NormErrorRepresentation} yield
		\begin{align*}
			\IntNormRB{u}{K}^2 - \IntNormFEM{u}{K}^2 &\preceq e^{-2C^*r}\sum\limits_{k=1}^N\left(\int_{0}^{1}\frac{t^{3-2s}\lambda_k^4}{1+t^2\lambda_k^2}u_k^2\,dt + \int_{1}^{\infty}\frac{t^{-2s-1}}{1+t^2\lambda_k^2}u_k^2\,dt\right) \\
			&\leq e^{-2C^*r}\sum\limits_{k=1}^Nu_k^2\left(\int_{0}^{1}\frac{t^{3-2s}\lambda_k^4}{t^2\lambda_k^2}\,dt + \int_{1}^{\infty}t^{-2s-1}\,dt\right)\\
			&\preceq e^{-2C^*r}\sum\limits_{k=1}^Nu_k^2\left(\lambda_k^2 - 1\right) \preceq e^{-2C^*r}\Norm{u}{1}^2.
		\end{align*}
	\end{proof}
	\subsection{Error of the reduced basis operator}\label{SecConvAna:SubSecOPError}
	The reduced basis interpolation norms provide an exponential decay in the error, granting good chances that similar results are valid with respect to the induced operators. Indeed, convergence of the operators is based on the results from Section \ref{SecConvAna:SubSecIntNormError}. The core of this paper is summarized in the following Theorem, relying on the notation
	\begin{align*}
		\Norm{v}{2}^2 := \sum\limits_{k=1}^N\lambda_k^{4}\scp{v}{\varphi_k}{0}^2, \rlap{\qquad$v\in\Vh$.}
	\end{align*}
	\begin{Th}[Exponential convergence of the reduced basis operator]\label{Th:Core}
		Let $u\in \Vh$ and $\Vr\subseteq \Vh$ a Zolotar\"ev space with $\sigma = [\lambda_L^2,\lambda_U^2]$ and $\delta = \nicefrac{\lambda_L^2}{\lambda_U^2}$. Then there exists a constant $C\in\R^+$, such that
		\begin{align}\label{SecConvAna:L2ConvOp}
			\|\InducedRBOp{H}(u) - \InducedFEMOp{H}u\|_{0} \preceq e^{-Cr}\Norm{u}{2}.
		\end{align}
		The constant $C$ only  depends on $\delta$ and satisfies
		\begin{align*}
			C(\delta) = \mathcal{O}\left(\frac{1}{\ln\left(\frac{1}{\delta}\right)}\right),\quad\text{as }\delta\to 0.
		\end{align*}
		Its precise value coincides with the constant $C^*$ from Remark \ref{Rem:AsympBehavC}. Moreover, if $s\in\left(0,\frac{1}{2}\right)$, the $2$-norm of $u$ in \eqref{SecConvAna:L2ConvOp} can be replaced by $\|u\|_{1}$.
	\end{Th}
	The rest of this section is dedicated to the proof of Theorem \ref{Th:Core} and therefore subject to the prescribed assumptions. All arising matrix-valued integrals are understood component-by-component.
	
	\begin{Lemma}\label{Lm:KsScalProd}
		There holds
		\begin{align*}
			\IntNormFEM{u}{K}^2 = \CoefVec{u}^TM\int_{0}^{\infty} t^{-2s-1}\left(M^{-1} - (M + t^2A)^{-1}\right)\,dt\,M\CoefVec{u}.
		\end{align*}
		Moreover, the induced scalar product $\InducedFEMScp{\cdot}{\cdot}{K}$ on $(\Vh,\IntNormFEM{\cdot}{K})$ satisfies
		\begin{align*}
			\scp{v}{w}{K^s} = \CoefVec{v}^TM\int_{0}^{\infty} t^{-2s-1}\left(M^{-1} - (M + t^2A)^{-1}\right)\,dt\,M\CoefVec{w}
		\end{align*}
		for all $v,w\in\Vh$.
	\end{Lemma}
	\begin{proof}
		We show the first equality. Due to \eqref{SecNormAppr:ShiftLapl}, there holds
		\begin{align*}
			\|u - \minimt{N}{t}\|_{0}^2 &=  \|\CoefVec{u} - (M + t^2A)^{-1}M\CoefVec{u}\|_M^2 \\
			&= \CoefVec{u}^TM\CoefVec{u} - 2\CoefVec{u}^TM(M + t^2A)^{-1}M\CoefVec{u} + \|(M + t^2A)^{-1}M\CoefVec{u}\|_M^2.
		\end{align*}
		Utilizing the identity $t^2A = (M+t^2A) - M$ yields
		\begin{align*}
			t^2\|\minimt{N}{t}\|_{1}^2 &= t^2\|(M + t^2A)^{-1}M\CoefVec{u}\|_{A}^2 \\
			&= \big((M + t^2A)^{-1}M\CoefVec{u}\big)^Tt^2A\big((M + t^2A)^{-1}M\CoefVec{u}\big) \\
			&= \CoefVec{u}^TM(M+t^2A)^{-1}M\CoefVec{u} - \|(M + t^2A)^{-1}M\CoefVec{u}\|_M^2.
		\end{align*}
		This reveals
		\begin{align*}
			K^2(t,u) &= \|u-\minimt{N}{t}\|_{0}^2 + t^2\|\minimt{N}{t}\|_{1}^2 \\
			&= \CoefVec{u}^TM\CoefVec{u} - \CoefVec{u}^TM(M+t^2A)^{-1}M\CoefVec{u} = \CoefVec{u}^TM(M^{-1} - (M + t^2A)^{-1})M\CoefVec{u}.
		\end{align*}
		Hence,
		\begin{align*}
			\IntNormFEM{u}{K}^2 &=  \int_{0}^{\infty} t^{-2s-1}K^2(t,u)\,dt \\
			&=\int_{0}^{\infty} t^{-2s-1}\CoefVec{u}^TM\left(M^{-1} - (M + t^2A)^{-1}\right)M\CoefVec{u}\,dt \\
			&= \CoefVec{u}^TM\int_{0}^{\infty} t^{-2s-1}\left(M^{-1} - (M + t^2A)^{-1}\right)\,dt\,M\CoefVec{u}.
		\end{align*}
	\end{proof}
	Throughout the rest of this paper, the identity matrix on $\R^{(r+1)\times(r+1)}$ is denoted by $I_r$, giving rise to the following claim.
	\begin{Lemma}\label{Lm:KsRbScalProd}
	There holds
		\begin{align*}
			\IntNormRB{u}{K}^2 = \CoefVec{u}^TMV_r\int_{0}^{\infty} t^{-2s-1}\left(I_r - (I_r + t^2A_r)^{-1}\right)\,dt\, V_r^TM\CoefVec{u}. 
		\end{align*}
	Moreover, the induced scalar product $\InducedRBScp{\cdot}{\cdot}{K}$ on $(\Vr,\IntNormRB{\cdot}{K})$ satisfies
	\begin{align}\label{SecConvAna:RBScp} 
		\scp{v_r}{w_r}{K_r^s} = \CoefVec{v_r}^TMV_r\int_{0}^{\infty} t^{-2s-1}\left(I_r - (I_r + t^2A_r)^{-1}\right)\,dt\, V_r^TM\CoefVec{w_r}
	\end{align}
	for all $v,w\in\Vh$. The matrix representation of $\InducedRBOp{K}$ is given by
	\begin{align*}
		\InducedRBMatOp{K} = MV_r\int_{0}^{\infty} t^{-2s-1}\left(I_r - (I_r + t^2A_r)^{-1}\right)\,dt\, V_r^TM.
	\end{align*}
	\end{Lemma}
	\begin{proof}
		In analogy to the proof of Lemma \ref{Lm:KsScalProd}, one shows that
		\begin{align*}
			K_r^2(t;u) = \CoefVec{u}^TMV_r(I_r - (I_r + t^2A_r)^{-1})V_r^TM\CoefVec{u}.
		\end{align*}
		Plugging into the integral representation of $\|u\|_{K_r^s}^2$ eventually proves the claim.
	\end{proof}
	\begin{Def}\label{Def:ExtendedSCP}
		For all $v\in\Vh$ and $w_r\in\Vr$ we define
		\begin{align*}
			\InducedRBScp{v}{w_r}{H}&:= \scp{\CoefVec{v}}{\InducedRBMatOp{H}\CoefVec{w_r}}{M}, \\
			\InducedRBScp{v}{w_r}{K} &:= \scp{\CoefVec{v}}{\InducedRBMatOp{K}\CoefVec{w_r}}{M}.
		\end{align*}  
	\end{Def}	
	Lemma \ref{Lm:KsScalProd} and \ref{Lm:KsRbScalProd} are fundamental to show pointwise convergence of the induced operator in the $\textrm{K}$-setting. We proceed in two steps.	
	\begin{Th}\label{Th:ExpConvForwOpPrelim}
		For all $w\in \Vh$ there holds
		\begin{align*}
			|\InducedRBScp{w}{u}{K} - \InducedFEMScp{w}{u}{K}| \leq \Norm{w}{0}\int_{0}^{\infty} t^{-2s-1} \sqrt{K_r^2(t;u) - K^2(t;u)}\,dt.
		\end{align*}
	\end{Th}
	\begin{proof}
		Due to Lemma \ref{Lm:KsScalProd} and \ref{Lm:KsRbScalProd}, $\InducedRBScp{w}{u}{K} - \InducedFEMScp{w}{u}{K}$ can be expressed as
		\begin{align*}
			\CoefVec{w}^TM\left(\int_{0}^{\infty} t^{-2s-1}\Big(V_r(I_r - (I_r + t^2A_r)^{-1})V_r^T - M^{-1} +(M+t^2A)^{-1}\Big)dt\right)M\CoefVec{u}.
		\end{align*}	
		One ascertains that the first term cancels out the third, i.e.,
		\begin{align*}
			\CoefVec{w}^TM\left(V_rV_r^T - M^{-1}\right)M\CoefVec{u} = \CoefVec{w}^T\left(MV_r\CoefvecVr{u} - M\CoefVec{u}\right) = 0.
		\end{align*}
		Computations of the remaining terms together with \eqref{SecNormAppr:LGS} and Remark \ref{Rem:ShiftedLaplaceRb} reveal
		\begin{align*}
			\CoefVec{w}^TM\left(-V_r(I_r + t^2A_r)^{-1}V_r^TM\CoefVec{u} + (M+t^2A)^{-1}M\CoefVec{u}\right) = \CoefVec{w}^TM\left(-\minimt{r}{t} + \minimt{N}{t}\right).
		\end{align*}
		Cauchy-Schwarz inequality leads to
		\begin{align*}
			|\InducedRBScp{w}{u}{K} - \InducedFEMScp{w}{u}{K}| &\leq \int_{0}^{\infty} t^{-2s-1} |\scp{w}{\minimt{N}{t}-\minimt{r}{t}}{0}|\,dt \\
		&\leq\Norm{w}{0}\int_{0}^{\infty} t^{-2s-1} \|\minimt{N}{t}- \minimt{r}{t}\|_{0}\,dt.
		\end{align*}
		There holds for all $t\in\R^+$
		\begin{align*}
			\|\minimt{N}{t} - \minimt{r}{t}\|_{0} &\leq \sqrt{\|\minimt{N}{t} - \minimt{r}{t}\|_{0}^2 + t^2\|\minimt{N}{t} - \minimt{r}{t}\|_{1}^2},
		\end{align*}
		which validates the conjecture by virtue of Corollary \ref{Cor:DiffKDiffInf}.
	\end{proof}
	\begin{Th}\label{Th:ExpConvForwardOperator}
			Let $C^*$ denote the constant from Remark \ref{Rem:AsympBehavC}. Then there holds
		\begin{align}\label{SecConvAna:ErrorSCP}
			\llap{$\forall w\in \Vh:$\quad}|\InducedRBScp{w}{u}{K} - \InducedFEMScp{w}{u}{K}| \preceq e^{-C^*r}\|w\|_{0} \Norm{u}{2}.
		\end{align}
		Moreover, if $s\in\left(0,\frac{1}{2}\right)$, the $2$-norm of $u$ in \eqref{SecConvAna:ErrorSCP} can be replaced by $\|u\|_{1}$.
	\end{Th}
	\begin{proof}
		We prove that for any $w\in \Vh$ and sufficiently small $\varepsilon>0$, satisfying $4s+2\varepsilon < 4$, there holds
		\begin{align*}
			|\InducedRBScp{w}{u}{K} - \InducedFEMScp{w}{u}{K}| \preceq e^{-C^*r} \|w\|_{0}\sqrt{\|u\|_{1}^2 + \sum\limits_{k=1}^N\lambda_k^{4s+2\varepsilon}u_k^2},
		\end{align*}
		which directly implies \eqref{SecConvAna:ErrorSCP}. Moreover, if $s<\frac{1}{2}$, we  can choose $\varepsilon<1-2s$ to verify the latter claim and conclude the proof. To this extent, let $\varepsilon\in(0,2-2s)$, if $s\geq\frac{1}{2}$, and $\varepsilon\in(0,1-2s)$ otherwise. Applying Theorem \ref{Th:ErrorKFunc2} followed by Cauchy-Schwarz inequality yields
		\begin{align*}
			\int_{0}^{1} t^{-2s-1}\sqrt{K_r^2(t;u) - K^2(t;u)}\,dt &\preceq \int_0^1t^{-\frac{1}{2}+\varepsilon}t^{-2s-\frac{1}{2}-\varepsilon}\sqrt{e^{-2C^*r}\sum\limits_{k = 1}^N\frac{t^4\lambda_k^4}{1+t^2\lambda_k^2}u_k^2}\,dt \\
			&\preceq e^{-C^*r}\sqrt{\int_0^1t^{-4s-1-2\varepsilon}\sum\limits_{k = 1}^N\frac{t^4\lambda_k^4}{1+t^2\lambda_k^2}u_k^2\,dt} \\
			&= e^{-C^*r}\sqrt{\sum\limits_{k = 1}^Nu_k^2\int_{0}^{1}\frac{t^{3-4s-2\varepsilon}\lambda_k^4}{1+t^2\lambda_k^2}\,dt}.
		\end{align*}
		Define $i := \min\{k\in\{1,...,N\}: \lambda_k \geq 1\}$ to observe
		\begin{align*}
			\sum\limits_{k = 1}^{i-1}u_k^2\int_{0}^{1}\frac{t^{3-4s-2\varepsilon}\lambda_k^4}{1+t^2\lambda_k^2}\,dt \leq \sum\limits_{k = 1}^{i-1} u_k^2\int_{0}^{\frac{1}{\lambda_k}}t^{3-4s-2\varepsilon}\lambda_k^4\,dt \preceq \sum\limits_{k = 1}^{i-1}u_k^2 \lambda_k^{4s+2\varepsilon}.
		\end{align*}
		Similarly, we obtain for the rest of the sum
				\begin{align*}
			\sum\limits_{k = i}^{N}u_k^2\int_{0}^{1}\frac{t^{3-4s-2\varepsilon}\lambda_k^4}{1+t^2\lambda_k^2}\,dt
			&\leq \sum\limits_{k = i}^{N}u_k^2\left( \int_{0}^{\frac{1}{\lambda_k}}t^{3-4s-2\varepsilon}\lambda_k^4\,dt + \int_{\frac{1}{\lambda_k}}^{1}\frac{t^{3-4s-2\varepsilon}\lambda_k^4}{t^2\lambda_k^2}\,dt\right)\\
			&\leq \sum\limits_{k = i}^{N}u_k^2\left( \frac{\lambda_k^{4s+2\varepsilon}}{4-4s-2\varepsilon} + \frac{\lambda_k^2-\lambda_k^{4s+2\varepsilon}}{2-4s-2\varepsilon}\right) \\
			&\preceq \|u\|_{1}^2 + \sum\limits_{k = i}^{N}u_k^2\lambda_k^{4s+2\varepsilon},
		\end{align*}
		such that
		\begin{align*}
			\int_{0}^{1} t^{-2s-1}\sqrt{K_r^2(t;u) - K^2(t;u)}\,dt\preceq e^{-C^*r}\sqrt{\|u\|_{1}^2 +\sum_{k=1}^Nu_k^2\lambda_k^{4s+2\varepsilon}}.
		\end{align*}
		On the interval $[1,\infty)$, we make use of Theorem \ref{Th:ErrorKFunc} to conclude for all $s\in(0,1)$
		\begin{align*}
			\int_{1}^{\infty} t^{-2s-1}\sqrt{K_r^2(t;u) - K^2(t;u)}\,dt &\preceq \int_1^\infty t^{-\frac{1}{2}-\varepsilon}t^{-2s-\frac{1}{2}+\varepsilon}\sqrt{e^{-2C^*r}\sum\limits_{k = 1}^N\frac{u_k^2}{1+t^2\lambda_k^2}}\,dt \\
			&\preceq e^{-C^*r}\sqrt{\sum\limits_{k = 1}^N u_k^2\int_1^\infty\frac{t^{-4s-1+2\varepsilon}}{1+t^2\lambda_k^2}\,dt} \\
			&\leq e^{-C^*r}\sqrt{\sum\limits_{k = 1}^N u_k^2\int_1^\infty \frac{t^{-4s-1+2\varepsilon}}{t^2\lambda_k^2}\,dt} \\
			&\preceq e^{-C^*r}\|u\|_{0}.
		\end{align*}
		Adding up the integrals in combination with Theorem \ref{Th:ExpConvForwOpPrelim} proves the claim.	
	\end{proof}
	As shown in the subsequent, equality of the scalar products is also valid with respect to Definition \ref{Def:ExtendedSCP}, providing the crucial link to transfer the error analysis form the $\textrm{K}$-setting to the spectral setting. 
	\begin{Lemma}\label{Lm:ScpEquality}
			For all $v\in\Vh$ and $w_r\in\Vr$ there holds
			\begin{align}\label{SecConvAna:ExtendedScpEquiv}
				\InducedRBScp{v}{w_r}{H} = C_s^2\InducedRBScp{v}{w_r}{K}. 
			\end{align}
	\end{Lemma}
	\begin{proof}
		If $v\in\Vr$, then \eqref{SecConvAna:ExtendedScpEquiv} holds due to Corollary \ref{Cor:ScpEquiv}. Let now $\Pi_r:\Vh\rightarrow\Vr$ denote the $0$-orthogonal projection on $\Vr$, such that
		\begin{align*}
			\llap{$\forall w_r\in\Vr:$\quad}\scp{\Pi_rv}{w_r}{0} = \scp{v}{w_r}{0}
		\end{align*}
		for each $v\in\Vh$. By virtue of the identity $\CoefvecVr{(\Pi_rv)}^T = \CoefVec{v}^TMV_r$, equation \eqref{SecOpAppr:ScpRep} and \eqref{SecConvAna:RBScp}, there holds for any arbitrary $v\in\Vh$ and $w_r\in\Vr$
		\begin{align*}
			\InducedRBScp{v}{w_r}{H} = \InducedRBScp{\Pi_rv}{w_r}{H} = C_s^2\InducedRBScp{\Pi_rv}{w_r}{K} = C_s^2\InducedRBScp{v}{w_r}{K}.
		\end{align*}
	\end{proof}
	We are eventually able to conduct the proof of Theorem \ref{Th:Core}.
	\begin{proof}[Proof of Theorem \ref{Th:Core}]
		Due to Corollary \ref{Cor:ScpEquiv} and Lemma \ref{Lm:ScpEquality}, there holds
		\begin{align*}
			\llap{$\forall w\in \Vh:$\quad}|\InducedRBScp{w}{u}{H} - \InducedFEMScp{w}{u}{H}|\leq |\InducedRBScp{w}{u}{K} - \InducedFEMScp{w}{u}{K}|.
		\end{align*}
		Hence,
		\begin{align*}
			\|\InducedRBOp{H}(u) - \InducedFEMOp{H}u\|_{0} &= \sup\limits_{w\in \Vh\setminus\{0\}}\frac{|\scp{w}{\InducedRBOp{H}(u) - \InducedFEMOp{H}u}{0}|}{\|w\|_{0}} \\
			&= \sup\limits_{w\in\Vh\setminus\{0\}}\frac{|\InducedRBScp{w}{u}{H} - \InducedFEMScp{w}{u}{H}|}{\|w\|_{0}} \\
			&\leq \sup\limits_{w\in\Vh\setminus\{0\}}\frac{|\InducedRBScp{w}{u}{K} - \InducedFEMScp{w}{u}{K}|}{\|w\|_{0}} \preceq e^{-C^*r}\Norm{u}{2},
		\end{align*}
		where the last inequality relies on \eqref{SecConvAna:ErrorSCP}. The fact that $\Norm{u}{2}$ reduces to $\Norm{u}{1}$ as $s<\frac{1}{2}$ follows directly from the latter claim of Theorem \ref{Th:ExpConvForwardOperator}.
	\end{proof}
	\begin{Rem}\label{Rem:EigenFuncContr}
		Results from Theorem \ref{Th:NormError} and \ref{Th:Core} can be ameliorated in a sense that the lower and upper bound, $\lambda_L^2$ and $\lambda_U^2$, solely have to be chosen with respect to that minimal and maximal eigenvalue, whose corresponding eigenfunction nontrivially contributes to the linear combination of the argument $u$. This leads to improvements of the constant $C^*$ and hence to faster convergence.
	\end{Rem}
	
	\section{Numerical examples}\label{SecNumExp}
	In the subsequent, we present several numerical examples in order to validate results from Section \ref{Sec:ConvAna}. To make matters precise, let $\Omega\subseteq\R^2$ be an open, bounded domain with Lipschitz boundary. We define $\Vnull := (L_2(\Omega),\|\cdot\|_{L_2})$ and $\Veins := (H_0^1(\Omega),\|\nabla\cdot\|_{L_2})$ to consider its arising interpolation space $\TheIntSpace$. The induced operator $\InducedOp{H}$ of $\IntNorm{\cdot}{H}$ coincides with the spectral fractional Laplacian subject to homogeneous Dirichlet boundary conditions, i.e.,
	\begin{align*}
		\llap{$\forall u\in\TheIntSpace:$\quad}\InducedOp{H}u = (-\Delta)^su := \sum_{k=1}^\infty \lambda_k^{2s}\scp{u}{\varphi_k}{L_2}\varphi_k,
	\end{align*}
	where $(\varphi_k,\lambda_k^2)_{k=1}^{\infty}\subseteq H_0^1(\Omega)\times\R^+$ denote the $L_2$-orthonormal eigenfunctions and eigenvalues of $\Vnull$ and $\Veins$. This setting is utilized to adequately study and analyze the performance of our method. All tests were implemented within the NGS-Py interface of the open source finite element packages NETGEN and NGSolve\footnote{www.ngsolve.org}, see \cite{Netgen} and \cite{NGSolve}. Computations of the Zolotar\"ev points are performed by means of the special function library from \texttt{Scipy}\footnote{https://docs.scipy.org/doc/scipy/reference/special.html}. 
	
	\begin{Ex}
	Consider the unit square $\Omega = (0,1)^2$ and a finite element space $\Vh\subseteq H_0^1(\Omega)$ of polynomial order $p  = 3$ on a quasi-uniform, triangular mesh $\mathcal{T}_h$ with mesh size $h = 0.08$, together with its arising eigenbasis $(\varphi_k)_{k=1}^N$ and $N = 1762$. We set
	\begin{align*}
		u = \sum\limits_{k=1}^n c_i\varphi_k
	\end{align*}
	for some randomly chosen coefficients $c_i\in(-1,1)$ with $n = 300$, such that the reduced basis norm is exact for $r\geq 299$. On $\Omega$, the exact eigenvalues $(\nu_k^2)_{k=1}^\infty$ of $-\Delta$ are given in closed form in terms of
	\begin{align*}
		\nu_k^2 = \nu_{i,j}^2 = \pi^2(i^2+j^2).
	\end{align*}
	Consistent with Remark \ref{Rem:EigenFuncContr}, we set $\lambda_L^2 := \nu_{1}^2\leq\lambda_1^2$ and utilize the power method to obtain an upper bound $\lambda_U^2 := 4200 \geq \lambda_{n}^2$. In accordance with Theorem \ref{Th:NormError}, Figure \ref{Fig:NormConvergence} affirms the exponential decay of the error $E_{u}^{Norm}(s,r) := \IntNormRB{u}{H}^2 - \IntNormFEM{u}{H}^2$ in $r$. Furthermore, we observe, as indicated in Section \ref{Sec:ConvAna}, that the error increases as $s$ is augmented.
	\begin{figure}[!h]
			\includegraphics[width=0.55\textwidth]{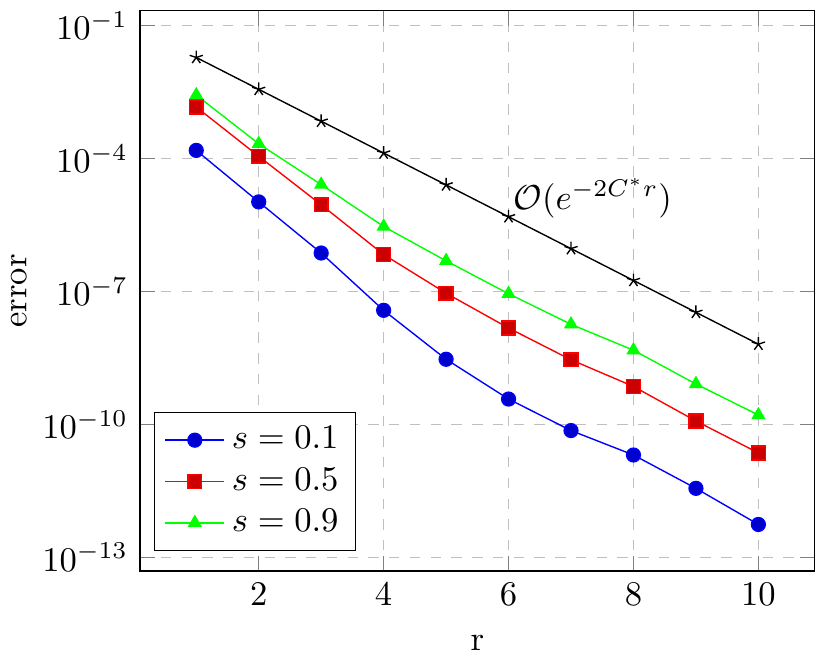}
			\caption{Error $E_{u}^{Norm}(s,r)$ of the reduced basis interpolation norm of a randomly chosen $u\in\Span\{\varphi_1,...,\varphi_n\}$ for three different values of $s$, $\sigma = [2\pi^2,4200]$, and $2C^*\approx 1.65$.}\label{Fig:NormConvergence}
	\end{figure}
	\end{Ex}
	
	\begin{Ex}
	Consider the unit circle $\Omega = \{\textbf{x}\in\R^2: \|\textbf{x}\|<1\}$, where $\|\cdot\|$ denotes the Euclidean norm, with corresponding finite element space $\Vh$ of polynomial order $p = 3$ on a quasi-uniform, triangular mesh $\mathcal{T}_h$, $h = 0.059$, with $N = 10270$. Let now $u_h := \Pi_h u\in \Vh$ denote the $L_2$-orthogonal projection of 
	\begin{align*}
		u(\textbf{x}) = u(x,y) = (1-\|\textbf{x}\|)y^2\sin(\|\textbf{x}\|)\in H_0^1(\Omega)
	\end{align*}
	onto $\Vh$. Having no further information which eigenfunctions contribute to the linear combination of $u_h$, we set $\lambda_L^2 := 1\leq \lambda_1^2$ and $\lambda_U^2 :=  \tilde{\lambda}_N^2+1\approx 1.53\cdot 10^5$, where $\tilde{\lambda}_N^2$ denotes a numerical approximation of $\lambda_N^2$ obtained by power iteration. The error of the reduced basis operator, $E_{u_h}^{Op}(s,r) := \|\InducedRBOp{H}(u_h) - \InducedFEMOp{H}u_h\|_{L_2}$, from Theorem \ref{Th:Core} is examined in Figure \ref{Fig:OperatorConvergence}. Here, the exact operator action $\InducedFEMOp{H}u_h$ has been replaced by $\Op_{H_{r^*}^s}(u_h)$, where $r^*\in\N$ is taken large enough to neglect the arising inaccuracy. The observed convergence rate is slightly better than $\mathcal{O}\left(e^{-C^*r}\right)$, where $C^*\approx 0.37$.
	
	\begin{figure}[ht]
		\begin{minipage}[t]{0.46\linewidth}
			\centering
			\includegraphics[width=\textwidth]{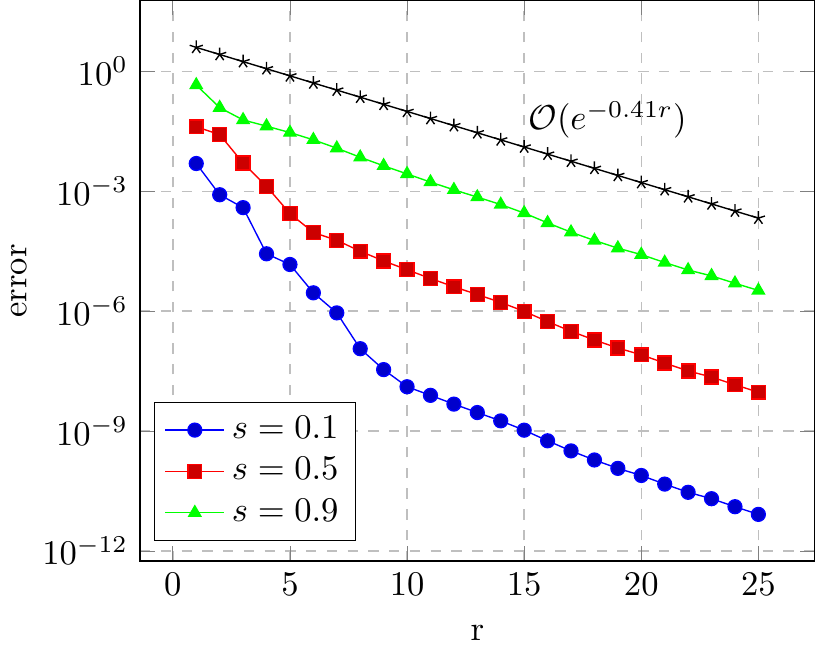}
			\captionsetup{width = \linewidth}
			\caption{$L_2$-error $E_{u_h}^{Op}(s,r)$ of the reduced basis operator for three different choices of $s$.}
			\label{Fig:OperatorConvergence}
		\end{minipage}
		\hspace{0.1cm}
		\begin{minipage}[t]{0.46\linewidth}
			\centering
			\includegraphics[width=\textwidth]{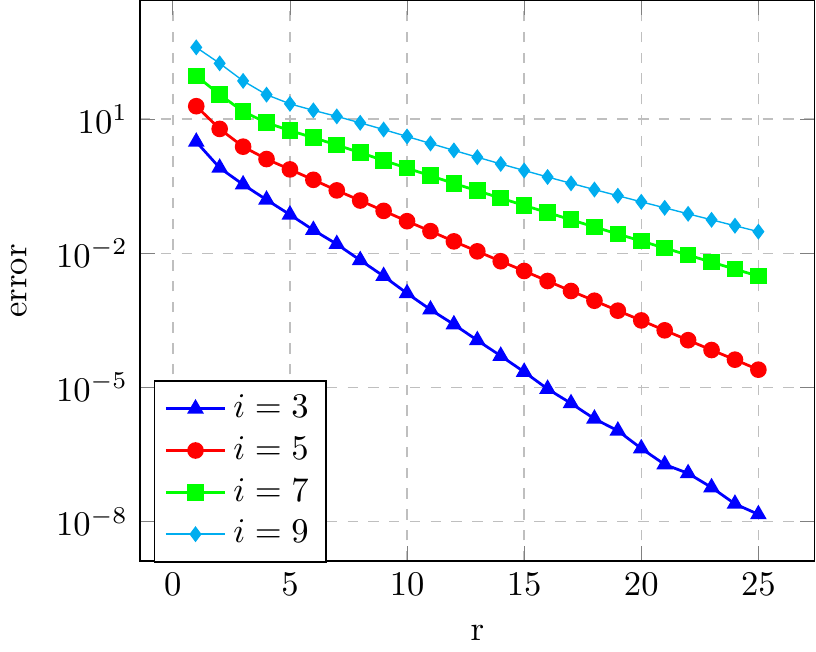}
			\captionsetup{width = \linewidth}
			\caption{Impact of mesh parameters $h = 2^{-i}$ on the error $E_{u_h^{rand}}^{Op}(0.9,r)$.}
			\label{Fig:Cond}
		\end{minipage}
	\end{figure}
	
	It is evident that the performance of our method relies on the condition of the problem and thus on the mesh parameter $h$. The exponential convergence property of $\InducedRBOp{H}$ for $h = 2^{-i}$, $i = 3,5,7,9$, on $\Omega$ with $p = 1$ and $s = 0.9$ is shown in Figure \ref{Fig:Cond} with respect to a randomly chosen $u_h^{rand}\in\Vh$. The rate of convergence deteriorates as $h\to0$. The error approximately behaves like $\mathcal{O}\left(e^{-0.77r}\right)$ for $h = 2^{-3}$ and $\mathcal{O}\left(e^{-0.31r}\right)$ for $h = 2^{-9}$.
	\end{Ex}
	
	\section{Appendix}
	\begin{proof}[Proof of Theorem \ref{Th:Trace}]
		The proof follows the outline of \cite[Proposition 2.1]{RadialExtremalSol} and \cite[Lemma 2.2]{ConConvEllipProb}. Let $v\in\Vbold$ with
		\begin{align*}
			v(y) = \sum\limits_{k=1}^\infty v_k(y)\varphi_k,
		\end{align*}
		and $v_k(y) = \scp{v(y)}{\varphi_k}{0}$. Then there holds
		\begin{align}\label{SecAppendix:MinimalExtNorm}
			\begin{split}
			\|v\|_{\mathbb{V}(\Vnull,\Veins; y^\alpha)}^2 &= \sum\limits_{k=1}^\infty\int_{\R^+}y^\alpha\left(\lambda_k^2 v_k(y)^2 + v_k'(y)^2\right)\,dy \\
			&\geq \sum\limits_{k=1}^\infty v_k(0)^2\min\limits_{\substack{\psi_k\in H^1(\R^+) \\ \psi_k(0) = 1}}\int_{\R^+}y^\alpha\left(\lambda_k^2 \psi_k(y)^2 + \psi_k'(y)^2\right)\,dy.
			\end{split}
		\end{align}
		It is well-known that the minimizer $\psi_k$ coincides with the solution of a Bessel-type differential equation
		\begin{align*}
			\begin{cases}
				\psi_k'' + \frac{\alpha}{y}\psi_k' - \lambda_k^2\psi_k = 0, &\text{in }\R^+, \\
				\lim\limits_{y\to \infty} \psi_k(y) = 0, \\
				\psi_k(0) = 1,
			\end{cases}
		\end{align*}
		which admits the following representation
		\begin{align*}
			\psi_k(y) = c_ky^s\textbf{K}_s(\lambda_ky).
		\end{align*} 
		Here, $c_k$ denotes a constant and $\textbf{K}_s$ the modified Bessel function of second kind, see \cite{HandbookOfMathFunc}. The constant is chosen in a way, such that $\psi_k(0) = 1$, which allows us to write 
		\begin{align}\label{SecInt:Psi}
			\psi_k(y) = \psi(\lambda_ky)
		\end{align}
		for a suitable function $\psi\in H^1(\R^+)$ with $\psi(0) = 1$. For all $k\in\N$ the value of the minimum in \eqref{SecAppendix:MinimalExtNorm} is given by
		\begin{align*}
			\int_{\R^+}y^{1-2s}\lambda_k^2\left(\psi(\lambda_ky)^2 + \psi'(\lambda_ky)^2\right)dy = \lambda_k^{2s}\int_{0}^\infty t^{1-2s}\left(\psi(t)^2 + \psi'(t)^2\right)dt.
		\end{align*}
		Integration by parts and incorporating the asymptotic behaviour of $\textbf{K}_\textsc{s}$ yields
		\begin{align*}
		\int_{0}^\infty t^{1-2s}\left(\psi(t)^2 + \psi'(t)^2\right)\,dt = d_s,
		\end{align*}
		see e.g., \cite[Remark 2.3]{ConConvEllipProb}. Thus,
		\begin{align*}
			\|v\|_{\Vbold}^2 \geq d_s\sum\limits_{k = 1}^\infty \lambda_k^{2s}v_k(0)^2 = d_s \IntNorm{v(0)}{H}^2.
		\end{align*}
		We conclude that $\Trace$ defines a linear, bounded operator that satisfies \eqref{SecInt:TraceIneq}, such that its range is contained in $\IntSpace{H}$. To evidence surjectivity, we observe that for each $u\in\IntSpace{H}$ the function
		\begin{align*}
			\mathscr{U}(y) := \sum\limits_{k=1}^\infty u_k\varphi_k\psi(\lambda_ky)
		\end{align*}
		is contained in $\Vbold$ and satisfies $\Trace\mathscr{U} = u$.
	\end{proof}

	\begin{proof}[Proof of Theorem \ref{Th:NormEquiv}]
		Let $u\in\IntSpace{E}$ and consider the ansatz function
		\begin{align*}
			\mathscr{U}(y) = \sum\limits_{k=1}^\infty u_k\varphi_k\psi(\lambda_ky),
		\end{align*}
		with $\psi$ as in \eqref{SecInt:Psi}. The proof of Theorem \ref{Th:Trace} validates that $\mathscr{U}\in\Vbold$ and 
		\begin{align*}
			\|\mathscr{U}\|_{\Vbold} = \sqrt{d_s}\IntNorm{u}{H}.
		\end{align*}
		Moreover, $\mathscr{U}$ satisfies \eqref{SecInt:HarmonicExtPDE}, by construction of $\psi$. What follows is that $\mathscr{U}$ is the $\alpha$-harmonic extension of $u$. Due to its minimization property from Lemma \ref{Lm:Minimizer}, we further have
		\begin{align*}
			\IntNorm{u}{E} = \|\mathscr{U}\|_{\Vbold} = \sqrt{d_s}\|u\|_{{\textrm{H}^s}(\Vnull,\Veins)}.
		\end{align*}
		To evidence the second equality, we refer to \cite[Theorem A.2]{Bramble}.
	\end{proof}

	\section*{Acknowledgements}
	One of the authors, Tobias Danczul, has been funded by the Austrian Science Fund (FWF) projects F 65 and W1245.

	\bibliographystyle{amsplain}
	\bibliography{Bibliography}

\end{document}